%% file: main.tex
\documentclass{article}
\usepackage[margin=1in]{geometry}
\usepackage{cite}
\usepackage[utf8]{inputenc}
\usepackage{amsmath,amssymb} %
\usepackage{xcolor}
\usepackage{float}
\usepackage{algorithm2e}
\usepackage{algorithmic}
\usepackage{subcaption}
\usepackage{csquotes}
\usepackage{url}
\usepackage{array}
\usepackage[thinlines]{easytable}
\usepackage{appendix}
\usepackage{ulem}
\usepackage{hyperref} 
\usepackage{booktabs}
\newcolumntype{L}{>{\centering\arraybackslash}m{4cm}}
\usepackage{lipsum}
\usepackage{graphicx}
\usepackage{breqn}
\graphicspath{{./Figures/}{../Figures}}
\RestyleAlgo{ruled}

\title{A Nonlinear Delay Model for Metabolic Oscillations in Yeast Cells}

\author{
    \textbf{Max M. Chumley}\\
    chumleym@msu.edu\\
    Michigan State University\\
    \and
    \textbf{Firas A.~Khasawneh$^*$}\\
    khasawn3@msu.edu\\
    Michigan State University\\
    \and
    \textbf{Andreas Otto}\\
    andreas.otto@iwu.fraunhofer.de\\
    Fraunhofer Institute\\ for Machine Tools\\ and Forming Technology IWU
    \and
    \textbf{Tomas Gedeon}\\
    gedeon@math.montana.edu\\
    Montana State University
}

\date{\today}

\begin{document} 

\maketitle
*Address all correspondence to this author.
\input{./Sections/abstract}
\input{./Sections/intro}

\input{./Sections/theory}

\input{./Sections/methods}
\input{./Sections/results}

\input{./Sections/conclusion}

\section{Acknowledgements}
This material is based upon work supported by the Air Force Office of Scientific Research under award number FA9550-22-1-0007.
Work of T.G was partially supported by NSF grant DMS-1951510.

\bibliographystyle{ieeetr}

\bibliography{Sections/bibliography}

\end{document}

%% file: Sections/abstract.tex
\section{Abstract}

We introduce two time-delay models of metabolic oscillations in yeast cells.  Our model tests a hypothesis that the oscillations occur as multiple pathways share a limited resource which we equate to the number of available ribosomes. We initially explore a single-protein model with a constraint equation governing the total resource available to the cell. The model is then extended to include three proteins that share a resource pool. Three approaches are considered at constant delay to numerically detect oscillations. First, we use a spectral element method to approximate the system as a discrete map and evaluate the stability of the linearized system about its equilibria by examining its eigenvalues. For the second method, we plot amplitudes of the simulation trajectories in 2D projections of the parameter space. We use a history function that is consistent with published experimental results to obtain metabolic oscillations. Finally, the spectral element method is used to convert the system to a boundary value problem whose solutions correspond to approximate periodic solutions of the system. Our results show that certain combinations of total resource available and the time delay, lead to oscillations. We observe that an oscillation region in the parameter space is between regions admitting steady states that correspond to zero and constant production. Similar behavior is found with the three-protein model where all proteins require the same production time. However, a shift in the protein production rates peaks occurs for low available resource suggesting that our model captures the shared resource pool dynamics.

%% file: Sections/intro.tex
\section{Introduction}

Cellular processes often exhibit non-trivial temporal dynamics in the absence of the external stimulus. Most common is the cell division cycle. However, as observed more than 50 years ago~\cite{Kuenzi1969}, yeast populations in low growth conditions exhibit metabolic cycling (MC)~\cite{Tu2005,Botstein2011} also known as respiratory cycling~\cite{Boczko}. While traditionally described as a result of carbon limitation, limitations by other essential nutrients like phosphate~\cite{Silverman} or ammonium, ethanol, phosphate, glucose, and sulfur~\cite{Brauer2008,Slavov2011} can lead to MC also known as metabolic oscillations. Under the growth conditions commonly used in this system, the population doubling time and thus the length of the cell division cycle is about $8$ h, and the metabolic oscillations have period $40-44$ min~\cite{Klavetz2004}.

The oscillations were first observed as periodic oscillations in the oxygen consumption of continuous, glucose-limited cultures growing in a chemostat, but were later also observed in batch cultures~\cite{Jules2005}. The MC has two distinct phases: low oxygen consumption (LOC) phase when dissolved oxygen in the medium is high and high oxygen consumption phase (HOC) when the oxygen in the medium drops to low levels~\cite{Botstein2011,Tu2005,Klavetz2004}.
Using experimental techniques ranging from micorarray analysis~\cite{Klavetz2004,Tu2005} to short-life luciferase fluorescent reporters~\cite{Boczko}, researchers were able to assign transcription of particular genes to these phases. During the LOC phase the yeast culture performs oxidative metabolism focused on amino-acid and ribosome synthesis, while during the HOC phase reductive reactions including DNA replication and proteosome related reactions occur~\cite{Tu2005}. Cellular metabolism during the reductive HOC phase seems to be devoted to the production of acetyl-CoA, preparing cells for the upcoming oxidative phase, during which metabolism shifts to respiration as accumulated acetyl-CoA units for ATP production via the TCA cycle and the electron transport chain~\cite{Tu2005}.

This compartmentalization of cellular processes in time is thought to be related to help assembly of macro-molecular complexes from units that, at low growth rates, are expressed at very low levels. Expressing them at the same time helps ensure timely synthesis and avoids waste of limited resources~\cite{Botstein2011}.

There were many different hypotheses centered on chemical signals that may mediate metabolic synchrony. 
In particular, Murray et. al.~\cite{Murray2003} proposed acetaldehyde and sulphate, Henson~\cite{Henson2004} and Sohn and Kuriyama~\cite{Sohn2001} hydrogen sulphide, while Adams et. al.~\cite{Adams2003}found that Gts1 protein plays a key stabilizing role. Finally, Muller at. al.~\cite{Muller2003} suggest a signalling agent, cAMP, plays a major role in mediating the integration of energy metabolism and cell cycle progression. 

Several mathematical models that do not specify the synchronizing chemical agent, but explore a general idea that cells in one phase of a cell cycle can slow down, or speed up progression of other cells through a different phase, have been suggested~\cite{boczko2010,stowers2011}. Finally, paper~\cite{young2018} explores synchronization which is a result of criticality of necessary cellular resources combined with the engagement of a cell cycle checkpoint, when these resources dip below the required level.

In this paper we explore the hypothesis that oscillations may arise spontaneously when several cellular processes share a limited resource. This does not explain why the processes separate into oxidative and reductive phase but argues that compartmentalization in time may help utilize limited resources more efficiently. Ribosomes are essential cellular resources as they produce enzymes used in all metabolic processes as well as all other proteins including those used to assemble ribosomes themselves. Yeast ribosomes are large molecular machines consisting of $79$ proteins~\cite{Woolford2013} and therefore they require substantial investment of cellular resources. This is reflected in the observation that the ratio of ribosomal proteins to all proteins scales linearly with cell growth rate across metabolic conditions~\cite{scott2014}. For this reason in our model we equate the limited shared cellular resource to the number of available ribosomes. 

In many biological systems, time delays are often incorporated into the models because many of these processes have nontrivial time spans that dictate the overall system behavior \cite{rihan2018applications, pell2023emergence,fowler1981approximate,gulbudak2021delay}. For this reason, the time delay framework is ideal for modeling a metabolic process where the protein production times can take upwards of 40 minutes. We aim to model the protein synthesis process in yeast cells using time delays and explore under what conditions oscillations are present in the responses. This paper is structured as follows. In Section~\ref{sec:theory_1ptn} we introduce a single protein model and extract theoretical results such as the fixed points and its linear stability behavior. Section~\ref{sec:theory_3ptn} presents the three protein extension to the single protein model and the equilibrium conditions are derived along with the system linearization. We then show the numerical methods that are utilized on the models in Section~\ref{sec:methods} where we describe the spectral element linear stability method, response feature analysis of system simulations under low growth conditions, and boundary value problem computation of periodic solutions to the nonlinear systems from simulation data. Results for the single protein system are then presented in Section~\ref{sec:1ptnresults} where the numerical methods are applied and the stability of the system is characterized in a subset of the overall parameter space. We then apply the same methods to the three protein system in Section~\ref{sec:3ptnresults}. Finally, we give concluding remarks in Section~\ref{sec:conclusion}.

%% file: Sections/theory.tex
\section{Theory - Single Protein}\label{sec:theory_1ptn}

\subsection{Model Derivation}
Both transcription and translation involve processing molecules (RNAP, ribosomes)
that are sequestered during the time of processing. These processing molecules constitute cellular resources that need to be shared by all necessary protein production processes. We will concentrate here on ribosomes, as their concentration is known to be tightly correlated with the microbial growth rate~\cite{scott2014}. The rate of production of protein $p(t)$ is proportional to the rate of initiation $\mu$ at some time $t-\tau(t)$ in the past when the processing started
\begin{equation} \label{eq:protein}
\dot{p}(t) = B \mu(t-\tau)  - D p(t),
\end{equation}
with the maximal growth rate $B$ and the decay rate $D$. The rate of initiation $\mu(t)$ is a product of the activator (which we assume for simplicity is the protein $p$ itself) and the ribosome $R$:
\begin{equation} \label{eq:initiation}
\mu(t) = f(p(t)) R(t), \\
\end{equation}
with Hill function 
\[f(t) = \frac{p^n(t)}{\kappa^n + p^n(t)}. \]
A suggestion for the sequestration equation based on \cite{TeranRomero2010} is given by 
\begin{equation} \label{eq:single_protein_resource}
R(t) = R_T - A \int_{t-\tau}^t \mu(s) ds, 
\end{equation}
where $R_T$ is the total resource (ribosomes) and the integral is the resource which is currently being sequestered to produce a protein. Differentiation of Eq.~\eqref{eq:single_protein_resource} leads to
\begin{equation} \label{eq:resource2}
\dot{R}(t) = A\left(\mu(t-\tau) - \mu(t)\right). 
\end{equation}
We note that this differentiation step is only valid for constant delays and if variable delays are used, Eq.~\eqref{eq:single_protein_resource} must be used for analysis.
Putting the equations together, we have the model
\begin{equation} \label{eq:single_protein_system}
\begin{split}
\dot{p}(t) &= B f(p(t-\tau)) R(t-\tau) - D p(t), \\
\dot{R}(t) &= A \left(f(p(t-\tau)) R(t-\tau) - f(p(t))) R(t)\right). \\
\end{split}
\end{equation}
The constant total resource $R_T$ is the sum of $R(t)$ and the integral over the history of $\mu(s)$ from $s=t-\tau(t)$ to $s=t$. This means that the initial functions for $p(\theta)$ and $R(\theta)$ with $\theta \in [-\tau(0),0]$ specify the value $R_T$.

\subsection{Equilibrium Points}

Equilibrium conditions of the system can be obtained by setting $p(t)=p(t-\tau)=p^*$, $R(t)=R(t-\tau)=R^*$, $\dot{p}(t)=0$, $\dot{R}(t)=0$ in Eq.~\eqref{eq:single_protein_system}. This process yields one equilibrium condition (Eq.~\eqref{eq:eq_cond1}) because the equation for $\dot{R}(t)$ in Eq.~\eqref{eq:single_protein_system} is satisfied for all constant $p$ and $R$. 

\begin{equation}\label{eq:eq_cond1}
	Dp^*=Bf(p^*)R^*
\end{equation}
The other equilibrium conditions are obtained from Eq.~\eqref{eq:single_protein_resource} by assuming that the integrand is constant when equilibrium has been reached. This yields Eq.~\eqref{eq:eq_cond2}.

\begin{equation}\label{eq:eq_cond2}
	R^*=\frac{R_T}{1+A\tau f(p^*)}
\end{equation}
We then substitute Eq.~\eqref{eq:eq_cond2} into Eq.~\eqref{eq:eq_cond1} to obtain,
\begin{equation}
	p^*=\frac{Bf(p^*)R_T}{D(1+A\tau f(p^*))}.
\end{equation}
Finally, inserting the definition of $f(p^*)$, we get a polynomial expression that must be satisfied for the equilibrium points as shown in Eq.~\eqref{eq:eq_cond3}.

\begin{equation}\label{eq:eq_cond3}
	(1+A \tau){p^*}^{n+1}-\frac{BR_T}{D}{p^*}^{n}+\kappa^n{p^*}=0
\end{equation}
Any solution to Eq.~\eqref{eq:eq_cond3} can then be plugged in to Eq.~\eqref{eq:eq_cond2} to obtain the equilibrium solutions. Note that this polynomial does not have analytical solutions for all values of $n$, but there is always one trivial equilibrium solution at $(p^*,R^*)=(0,R_T)$. The nontrivial equilibrium points are then obtained from solving the following system for ($p^*,R^*$).

\begin{subequations}
    \begin{equation}\label{eq:eq_poly}
        (1+A \tau){p^*}^{n}-\frac{BR_T}{D}{p^*}^{n-1}+\kappa^n=0,
    \end{equation}
        	
    \begin{equation}
        R^*=\frac{R_T}{1+A\tau f(p^*)}.
    \end{equation}
\end{subequations}
Once the trivial solution is removed from the conditions, the remaining polynomial in $p^*$ has either 0 or 2 positive real roots by Descartes' rule of signs for every $~n\in \mathbb{Z}^+$ assuming only positive parameters are chosen. In conclusion, there is at least 1 trivial root at $(0,R_T)$ and at most 3 equilibrium points including the trivial point where the other two points are the nontrivial equilibria. We choose to restrict the analysis to $n=2$ so that we can obtain analytical solutions for the equilibrium points. Solving Eq.~\eqref{eq:eq_cond3} for $n=2$ yielded three equilibrium points:

\begin{equation*}
	\begin{split}
        &\left(0,R_T\right),\\
        \\
        &\left(\frac{B R_T-\sqrt{B^2 R_T^2-4 D^2 \kappa ^2 (A \tau +1)}}{2 D(A  \tau +1)}, \frac{B R_T \sqrt{B^2 R_T^2-4 D^2 \kappa ^2 (A \tau +1)}-2 A D^2 \kappa ^2 \tau  (A \tau +1)-B^2 R_T^2}{B (A \tau +1) \left(\sqrt{B^2 R_T^2-4 D^2 \kappa ^2 (A \tau +1)}-B R_T\right)}\right), \\
        \\
        \text{and}
        \\
        &\left(\frac{B R_T+\sqrt{B^2 R_T^2-4 D^2 \kappa ^2 (A \tau +1)}}{2 D(A  \tau +1)}, \frac{B R_T \sqrt{B^2 R_T^2-4 D^2 \kappa ^2 (A \tau +1)}+2 A D^2 \kappa ^2 \tau  (A \tau +1)+B^2 R_T^2}{B (A \tau +1) \left(\sqrt{B^2 R_T^2-4 D^2 \kappa ^2 (A \tau +1)}+B R_T\right)}\right),
	\end{split}
	\label{eq:single_eq_pts}
\end{equation*}
where, $A, B, D, \tau, \kappa, \text{ and } R_T$ are system parameters. Note that the second equilibrium point has a protein production rate that is between the protein production rates of the other two equilibrium points. Therefore, we refer to this equilibrium point as the ``middle" equilibrium point and the point with the largest $p^*$ as the ``top" equilibrium point. By studying the stability of these three fixed points, the stability of the system can be characterized for certain parameters.

To understand the role that each equilibrium point plays in the system, we need to understand which equilibrium points are present in different regions of the parameter space. For this system, we can compute the saddle node bifurcation by studying the curve where the argument in the square root term becomes negative indicating that the top and middle equilibria are no longer real numbers. This curve forms a boundary that separates the region with three equilibrium points and the region with only the trivial equilibrium point. It can be computed analytically for this system and is defined in terms of $\tau$ and $R_T$ as:
\begin{equation}\label{eq:saddlenode}
	R_T<\sqrt{\frac{4D^2 \kappa^2(1+A\tau)}{B^2}}.
\end{equation}
This boundary given by the equality of Eq.~\eqref{eq:saddlenode} is plotted in the stability diagrams in Sec.~\ref{sec:1ptnresults} as a red curve with triangles. So any parameters that satisfy \eqref{eq:saddlenode} only have a single equilibrium at the trivial point, and if the parameters do not satisfy this inequality, the top and middle equilibrium points are valid equilibrium solutions along with the trivial point.

\subsection{System Linearization}\label{sec:single_lin}
In the analysis of nonlinear dynamical systems it is useful to study the associated linearized system about the fixed points. The Hartman-Grobman theorem states that near a hyperbolic equilibrium point, the linearized system exhibits the same behavior as the nonlinear system \cite{seydel2009practical}. We linearize \eqref{eq:single_protein_system} by computing the Jacobian matrices of the present and delayed states about an equilibrium point $\vec{q}=\begin{bmatrix}p^*~R^*\end{bmatrix}^T$. We start by defining two state space vectors, $\vec{x}$ and $\vec{x}_\tau$
where,

\begin{equation*}
	\begin{split}
	\vec{x} &= \begin{bmatrix}
		p(t)\\
		R(t)
		\end{bmatrix}=\begin{bmatrix}
			x_1\\
			x_2
			\end{bmatrix},\\ 
			\vec{x}_\tau &= \begin{bmatrix}
				p(t-\tau)\\
				R(t-\tau)
				\end{bmatrix}=\begin{bmatrix}
					x_{1 \tau}\\
					x_{2 \tau}
					\end{bmatrix}.
		\end{split}
\end{equation*}
The nonlinear delay system can then be written in the form,
\begin{equation}\label{eq:g_single}
	\dot{\vec{x}}(t)=\vec{g}(\vec{x}, \vec{x}_\tau),
\end{equation}
where,
$\vec{g}=\vec{g}_1(\vec{x})+\vec{g}_2(\vec{x}_\tau)$,~
$\vec{g}_1=\begin{bmatrix}
	-Dx_1\\
	-Af(x_1)x_2
\end{bmatrix}$,~
$\vec{g}_2=\begin{bmatrix}
	Bf(x_{1\tau})x_{2\tau}\\
	Af(x_{1\tau})x_{2\tau}
\end{bmatrix}$. In this form, the system is written as a sum of a nonlinear component as a function of $t$ and a nonlinear delay component as a function of $t-\tau$. We linearize each piece of $g$ by computing the Jacobian matrix of the vector functions.
\begin{equation}
	\textbf{G}_1=\frac{\partial \vec{g}_{1}}{\partial x}=\begin{bmatrix}
		\frac{\partial \vec{g}_{11}}{\partial x_1} & \frac{\partial \vec{g}_{11}}{\partial x_2}\\
		\frac{\partial \vec{g}_{12}}{\partial x_1} & \frac{\partial \vec{g}_{12}}{\partial x_2}\\
	\end{bmatrix}=\begin{bmatrix}
		-D & 0\\
		-A f'(x_1)x_2 & -A f(x_1)
	\end{bmatrix},
\end{equation}
and
\begin{equation}
	\textbf{G}_2=\frac{\partial \vec{g}_{2}}{\partial x_\tau}=\begin{bmatrix}
		\frac{\partial \vec{g}_{21}}{\partial x_{1\tau}} & \frac{\partial \vec{g}_{21}}{\partial x_{2\tau}}\\
		\frac{\partial \vec{g}_{22}}{\partial x_{1\tau}} & \frac{\partial \vec{g}_{22}}{\partial x_{2\tau}}\\
	\end{bmatrix}=\begin{bmatrix}
		B f'(x_{1\tau})x_{2\tau} & B f(x_{1\tau})\\
		Af'(x_{1\tau})x_{2\tau} & Af(x_{1\tau})
	\end{bmatrix},
\end{equation} where $f'(x_{1\tau})=\frac{n\kappa^nx_{1\tau}^{n-1}}{(\kappa^n+x_{1\tau}^n)^2}$. The linearized system about an equilibrium $q$ is then written as

\begin{equation}\label{eq:lin_sys1}
	\dot{\vec{x}}\approx\begin{bmatrix}
		-D & 0\\
		-A f'(x_1)x_2 & -A f(x_1)
	\end{bmatrix}\bigg|_q (\vec{x}-\vec{q})+\begin{bmatrix}
		B f'(x_{1\tau})x_{2\tau} & B f(x_{1\tau})\\
		Af'(x_{1\tau})x_{2\tau} & Af(x_{1\tau})
	\end{bmatrix}\bigg|_q(\vec{x}_\tau-\vec{q}).
\end{equation}
If a change of variables, $\vec{y} = \vec{x} - \vec{q}$ is implemented, with $\vec{y}_\tau = \vec{x}_\tau - \vec{q}$, Eq.~\eqref{eq:lin_sys1} simplifies to Eq.~\eqref{eq:lin_sys2} effectively moving the equilibrium point to the origin.
\begin{equation}\label{eq:lin_sys2}
	\dot{\vec{y}}\approx\begin{bmatrix}
		-D & 0\\
		-A f'(p^*)R^* & -A f(p^*)
	\end{bmatrix} \vec{y}+\begin{bmatrix}
		B f'(p^*)R^* & B f(p^*)\\
		Af'(p^*)R^* & Af(p^*)
	\end{bmatrix}\vec{y}_\tau.
\end{equation}
We can write Eq.~\eqref{eq:lin_sys2} in a simplified form as,
\begin{equation}\label{eq:lin_sys_simplified}
	\dot{\vec{y}}\approx\textbf{G}_1(q)\vec{y}+\textbf{G}_2(q)\vec{y}_\tau.
\end{equation}
We will use this linearized system to evaluate the stability of the equilibrium points of the system. Note that Eq.~\eqref{eq:lin_sys_simplified} was derived only from the DDE system Eq.~\eqref{eq:single_protein_system}. In addition, the constraint Eq.~\eqref{eq:single_protein_resource} must hold also for the perturbations. It is straightforward to directly compute the linearized system about the trivial equilibrium. We do this by inserting $(p^*,R^*)=(0,R_T)$ into \eqref{eq:lin_sys_simplified}, which leads to the elementary ODE system 
\begin{equation}\label{eq:lin_sys_trivial}
	\dot{\vec{y}}\approx\begin{bmatrix}
		-D & 0\\
		0 & 0
	\end{bmatrix} \vec{y},
\end{equation}
because all elements of $\textbf{G}_2$ become zero. This system has only two characteristic exponents that are directly obtained from the diagonal of $\textbf{G}_1$ as $-D$ and $0$. Since the original nonlinear DDE system is infinite dimensional, there are infinitely many other characteristic exponents, which all tend to $-\infty$ as $\vec{q} \to \begin{bmatrix}0~R_T\end{bmatrix}^T$. Moreover, the characteristic exponent equal to zero corresponds to perturbations of the resources $R(t)$, which changes the value of the overall resources $R_T$. However, such perturbations do not fulfill the additional constraint equation \eqref{eq:single_protein_resource}, which means that this eigenvalue corresponds to the eigenvector along a one dimensional family of  equilibria parameterized by the total resource $R_T$. As such this eigenvalue does not reflect the stability of the equilibrium within a phase space with $R_T$ fixed. As a result, the trivial equilibrium point is a locally stable node as long as the decay rate is positive $(D>0)$. We emphasize that using the Jacobian methods for linearization at this step is valid for constant delays. For state-dependent delays or time-dependent delays the system can be linearized using methods from \cite{gedeon2022operon}.

\section{Theory - Three Protein Model}\label{sec:theory_3ptn}
\subsection{Model}
We extend the single protein model in Eq.~\eqref{eq:single_protein_system} to incorporate production of three proteins with shared resource i.e. a shared ribosomal pool. If the resources are shared, we expect that oscillations may occur if there are not enough resources to produce all three proteins simultaneously. The extended model is shown in Eq.~\eqref{eq:three_protein_system}.

\begin{equation}
    \begin{split}
        \dot{p_1}(t)&=B_1 f(p_2(t-\tau_1))f(p_3(t-\tau_1))R(t-\tau_1) - D_1 p_1,\\
        \dot{p_2}(t)&=B_2 f(p_1(t-\tau_2))R(t-\tau_2) - D_2 p_2,\\
        \dot{p_3}(t)&=B_3 f(p_1(t-\tau_3))R(t-\tau_3) - D_3 p_3,\\
        \dot{R}(t)&=A(\mu_1(t-\tau_1)+\mu_2(t-\tau_2)+\mu_2(t-\tau_3)-\mu_1(t)-2\mu_2(t)),
    \end{split}
    \label{eq:three_protein_system}
\end{equation}
where $A$, $B_1$, $B_2$, $B_3$, $\tau_1$, $\tau_2$, $\tau_3$, $D_1$, $D_2$, $D_3$, are system parameters, $\mu_1(t)=f(p_2(t))f(p_3(t))R(t)$, and $\mu_2(t)=f(p_1(t))R(t)$ and $f(x)=\frac{x^n}{\kappa^n+x^n}$. This means that production of the first protein is activated when the other two protein production rates are nonzero and production of the second and third proteins is activated by $p_1$. Analogously to the single protein system, the total resource ($R_T$) is computed using,
\begin{equation}
    R_T = R(t) + A\left(\int_{t-\tau_1}^t f(p_2(s))f(p_3(s))R(s)ds + \int_{t-\tau_2}^t f(p_1(s))R(s)ds + \int_{t-\tau_3}^t f(p_1(s))R(s)ds\right). \label{eq:rt}
\end{equation}

\subsection{Equilibrium Points}
The equilibrium conditions are found by first setting $\dot{p_1}=\dot{p_2}=\dot{p_3}=0$ yielding the following conditions:
\begin{equation}
    \begin{split}
        D_1p_1^*&=B_1 f(p_2^*)f(p_3^*)R^*,\\
    	D_2p_2^*&=B_2f(p_1^*)R^*,\\
    	D_3p_3^*&=B_3f(p_1^*)R^*,
    \end{split}
     \label{eq:eq_cond}
\end{equation} 
where ($p_1^*$, $p_2^*$, $p_3^*$,$R^*$) is the equilibrium point. Similarly to the single protein system, the $\dot{R}$ expression in Eq.~\eqref{eq:three_protein_system} is always satisfied at equilibrium, but Eq.~\eqref{eq:rt} yields the fourth and final equilibrium condition:
\begin{equation}
	R_T = R^*+A\left[f(p_2^*)f(p_3^*)R^*\tau_1+f(p_1^*)R^*(\tau_2+\tau_3)\right]. \label{eq:rteq}
\end{equation}

This system has a trivial equilibrium at $p_1^*=p_2^*=p_3^* = 0, R^*=R_T$. For finding the other equilibria, we need to solve this system of equations. We see that the relations in Eq.~\eqref{eq:eq_cond} all depend directly on $R^*$. We eliminate $R^*$ by solving Eq.~\eqref{eq:rteq} for $R^*$,
\begin{equation}\label{eq:r_star}
	R^* = \frac{R_T}{1+A\left(f(p_2^*)f(p_3^*)\tau_1+f(p_1^*)(\tau_2+\tau_3)\right)},
\end{equation}
and substituting $R^*$ in the three Eqs.~\eqref{eq:eq_cond}, along with using the definition of $f(x)$. These steps result in three multivariate polynomial equilibrium equations shown in Eqs.~\eqref{eq:poly_eq_1},\eqref{eq:poly_eq_2} and \eqref{eq:poly_eq_3}.

\begin{dmath}\label{eq:poly_eq_1}
	D_1p_1^*(\kappa^n+p_1^{*n})(\kappa^n+p_2^{*n})(\kappa^n+p_3^{*n})+AD_1p_1^*p_2^{*n}p_3^{*n}(\kappa^n+p_1^{*n})\tau_1+AD_1p_1^{*(n+1)}(\kappa^n+p_2^{*n})(\kappa^n+p_3^{*n})(\tau_2+\tau_3) = B_1R_Tp_2^{*n}p_3^{*n}(\kappa^n+p_1^{*n}),
\end{dmath}
\begin{dmath}\label{eq:poly_eq_2}
	D_2p_2^*(\kappa^n+p_1^{*n})(\kappa^n+p_2^{*n})(\kappa^n+p_3^{*n})+AD_2p_2^{*(n+1)}p_3^{*n}(\kappa^n+p_1^{*n})\tau_1+AD_2p_1^{*n}p_2^*(\kappa^n+p_2^{*n})(\kappa^n+p_3^{*n})(\tau_2+\tau_3) = B_2R_Tp_1^{*n}(\kappa^n+p_2^{*n})(\kappa^n+p_3^{*n}),
\end{dmath}
\begin{dmath}\label{eq:poly_eq_3}
	D_3p_3^*(\kappa^n+p_1^{*n})(\kappa^n+p_2^{*n})(\kappa^n+p_3^{*n})+AD_3p_2^{*n}p_3^{*(n+1)}(\kappa^n+p_1^{*n})\tau_1+AD_3p_1^{*n}p_3^*(\kappa^n+p_2^{*n})(\kappa^n+p_3^{*n})(\tau_2+\tau_3) = B_3R_Tp_1^{*n}(\kappa^n+p_2^{*n})(\kappa^n+p_3^{*n}).
\end{dmath}
The solutions $(p_1^*, p_2^*, p_3^*)$ to these three equations correspond to the equilibrium point of the system and the resource equilibrium is obtained by substituting these values in to Eq.~\eqref{eq:r_star}. Due to the complexity of these equations, we solve them numerically using the variable precision (VPA) solver in \texttt{Matlab} \cite{vpasolve}. Details for how these equations were solved are outlined in Sec.~\ref{sec:3ptnresults}.

\subsection{Three Protein System Linearization}\label{sec:three_lin}
The three protein system was also linearized about its equilibrium points for stability analysis using the spectral element linear stability method described in section \ref{sec:spec_stability}. We will linearize the system about the equilibrium point, $\vec{q}=\begin{bmatrix}p_1^*~p_2^*~p_3^*~R^*\end{bmatrix}^T$ from the solution to Eqs.~\eqref{eq:poly_eq_1},\eqref{eq:poly_eq_2}, and \eqref{eq:poly_eq_3}. Similarly to the single protein model, we define state vectors for the current states and delayed states, but in this case, three delayed states are present due to the system having multiple time delays. The system states are,

\begin{align}
	\vec{x} &= \begin{bmatrix}
		p_1(t)\\
		p_2(t)\\
		p_3(t)\\
		R(t)
		\end{bmatrix}=\begin{bmatrix}
			x_1\\
			x_2\\
			x_3\\
			x_4
			\end{bmatrix},~
			\vec{x}_{\tau_i} = \begin{bmatrix}
				p_1(t-\tau_i)\\
				p_2(t-\tau_i)\\
				p_3(t-\tau_i)\\
				R(t-\tau_i)
				\end{bmatrix}=\begin{bmatrix}
					x_{1 \tau_i}\\
					x_{2 \tau_i}\\
					x_{3 \tau_i}\\
					x_{4 \tau_i}
					\end{bmatrix},
\end{align}
where $i \in [1,2,3]$ represents the system state at delay $\tau_i$. We then write the system as:
\begin{equation}\label{eq:g_triple}
	\dot{\vec{x}} = \vec{g}(\vec{x},\vec{x}_{\tau_1},\vec{x}_{\tau_2},\vec{x}_{\tau_3}),
\end{equation}
and separate $\vec{g}$ as a sum of terms only dependent on one of the system states.
\begin{equation*}
	\vec{g}(\vec{x},\vec{x}_{\tau_1},\vec{x}_{\tau_2},\vec{x}_{\tau_3}) = \vec{g}_1(\vec{x}) + \vec{g}_2(\vec{x}_{\tau_1}) + \vec{g}_3(\vec{x}_{\tau_2}) + \vec{g}_4(\vec{x}_{\tau_3}),
\end{equation*}
where,
\begin{align*}
	&\vec{g}_1(\vec{x}) = \begin{bmatrix}
		-D_1x_1\\
		-D_2x_2\\
		-D_3x_3\\
		-Af(x_1)x_4-2Af(x_2)x_4
	\end{bmatrix},\indent
	&\vec{g}_2(\vec{x}_{\tau_1}) =& \begin{bmatrix}
		B_1f(x_{2\tau})f(x_{3\tau})x_{4\tau}\\
		0\\
		0\\
		Af(x_{1\tau})x_{4\tau}
	\end{bmatrix},\\
	&\vec{g}_3(\vec{x}_{\tau_2}) = \begin{bmatrix}
		0\\
		B_2f(x_{1\tau})x_{4\tau}\\
		0\\
		Af(x_{2\tau})x_{4\tau}
	\end{bmatrix},\indent
	&\vec{g}_4(\vec{x}_{\tau_3}) =& \begin{bmatrix}
		0\\
		0\\
		B_3f(x_{1\tau})x_{4\tau}\\
		Af(x_{2\tau})x_{4\tau}
	\end{bmatrix},
\end{align*}
where $f(x)=\frac{x^n}{\kappa^n+x^n}$. Now, the system can be linearized about $q$ as,
\begin{equation}\label{eq:lin_x}
	\dot{\vec{x}}\approx \textbf{G}_1(\vec{q})(\vec{x}-\vec{q})+\sum_{i=2}^4\textbf{G}_i(q)(\vec{x}_{\tau_i} - \vec{q}),
\end{equation}
where $\textbf{G}_i$ is the Jacobian matrix of $\vec{g}_i(\vec{x}_{\tau_{i-1}})$. The Jacobian matrices for this system are analytically computed with the matrix of partial derivatives. For example, $\textbf{G}_1$ is computed as,
\begin{equation*}
	\textbf{G}_1 = \begin{bmatrix}
		\frac{\partial \vec{g}_{11}}{\partial x_1} &\frac{\partial \vec{g}_{11}}{\partial x_2} &\frac{\partial \vec{g}_{11}}{\partial x_3}& \frac{\partial \vec{g}_{11}}{\partial x_4}\\[6pt]
		\frac{\partial \vec{g}_{12}}{\partial x_1} &\frac{\partial \vec{g}_{12}}{\partial x_2} &\frac{\partial \vec{g}_{12}}{\partial x_3}& \frac{\partial \vec{g}_{12}}{\partial x_4}\\[6pt]
		\frac{\partial \vec{g}_{13}}{\partial x_1} &\frac{\partial \vec{g}_{13}}{\partial x_2} &\frac{\partial \vec{g}_{13}}{\partial x_3}& \frac{\partial \vec{g}_{13}}{\partial x_4}\\[6pt]
		\frac{\partial \vec{g}_{14}}{\partial x_1} &\frac{\partial \vec{g}_{14}}{\partial x_2} &\frac{\partial \vec{g}_{14}}{\partial x_3}& \frac{\partial \vec{g}_{14}}{\partial x_4}\\
	\end{bmatrix},
\end{equation*}
where $\frac{\partial g_{1j}}{\partial x_k}$ is the partial derivative of the $j$-th component of the $g_1$ vector with respect to $x_k$. A similar process is used for $\textbf{G}_2,~\textbf{G}_3~\text{and}~\textbf{G}_4$ with the main difference being that the derivatives are computed with respect to delayed states at the corresponding delay $\tau_i$. Carrying out this procedure yields the following expressions for the Jacobian matrices.

\begin{align*}
	\textbf{G}_1(\vec{q}) &= \begin{bmatrix}
		-D_1 & 0 & 0 & 0\\
		0 & -D_2 & 0 & 0\\
		0 & 0 & -D_3 & 0\\
		-2Af'(p_1^*)R^* & -Af'(p_2^*)f(p_3^*)R^* & -Af(p_2^*)f'(p_3^*)R^* & -2Af(p_1^*)-Af(p_2^*)f(p_3^*)
	\end{bmatrix},\\[6pt]
	\textbf{G}_2(\vec{q}) &= \begin{bmatrix}
		0 & B_1f'(p_2^*)f(p_3^*)R^* & B_1f(p_2^*)f'(p_3^*)R^* & B_1f(p_2^*)f(p_3^*)\\
		0 & 0 & 0 & 0\\
		0 & 0 & 0 & 0\\
		0 & Af'(p_2^*)f(p_3^*)R^* & Af(p_2^*)f'(p_3^*)R^* & Af(p_2^*)f(p_3^*)
	\end{bmatrix},\\[6pt]
	\textbf{G}_3(\vec{q}) &= \begin{bmatrix}
		0 & 0 & 0 & 0\\
		B_2f'(p_1^*)R^* & 0 & 0 & B_2f(p_1^*)\\
		0 & 0 & 0 & 0\\
		Af'(p_1^*)R^* & 0& 0 & Af(p_1^*)
	\end{bmatrix},\\[6pt]
	\textbf{G}_4(\vec{q}) &= \begin{bmatrix}
		0 & 0 & 0 & 0\\
		0 & 0 & 0 & 0\\
		B_3f'(p_1^*)R^* & 0 & 0 & B_3f(p_1^*)\\
		Af'(p_1^*)R^* & 0& 0 & Af(p_1^*)
	\end{bmatrix},
\end{align*}
where $f'(x)=\frac{n\kappa^nx^{n-1}}{(\kappa^n+x^n)^2}$. Finally, we introduce the change of variables $\vec{y}=\vec{x}-\vec{q} $ to Eq.~\eqref{eq:lin_x} resulting in the linearized system:
\begin{equation}\label{eq:3ptn_linear}
	\dot{y} \approx \textbf{G}_1(\vec{q})\vec{y}(t) + \textbf{G}_2(\vec{q})\vec{y}(t-\tau_1)+ \textbf{G}_3(\vec{q})\vec{y}(t-\tau_2) + \textbf{G}_4(\vec{q})\vec{y}(t-\tau_3). 
\end{equation} 

For the trivial equilibrium $(0,0,0,R_T)$, we have $f(0)=0$ and $f'(0)=0$, and the matrices $\textbf{G}_2$, $\textbf{G}_3$ and $\textbf{G}_4$ vanish. Thus, similarly to the single protein system, the linearization at the trivial equilibrium becomes a simple ODE system 
\begin{equation*}
    \dot{y}\approx\begin{bmatrix}
		-D_1 & 0 & 0 & 0\\
		0 & -D_2 & 0 & 0\\
		0 & 0 & -D_3 & 0\\
		0 & 0 & 0 & 0
	\end{bmatrix}\vec{y}.
\end{equation*}
This system has four eigenvalues $-D_1$, $-D_2$, $-D_3$, and zero. Again, the characteristic exponent equal to zero corresponds to family of equilibria parameterized by $R_T$ and can be discarded. We conclude that the trivial equilibrium point is locally stable for all positive decay rates ($D_1>0,~D_2>0,~D_3>0$).

%% file: Sections/methods.tex
\section{Methods} \label{sec:methods}

 \subsection{Spectral Element Approach - Linear Stability Analysis}\label{sec:spec_stability}

For analyzing the dynamic behavior of the system and identify regions of bistablity in parameter space, we study the linear stability of the equilibria. We use the spectral element method, which is an advanced numerical method for the stability analysis of DDE systems \cite{Khasawneh2013}. In particular, the linear variational systems Eq.~\eqref{eq:lin_sys_simplified} and Eq.~\eqref{eq:3ptn_linear} for perturbations around the equilibrium points are converted to a dynamic map, which describes the evolution of the system state $\vec{z}_{n-1}$ at time step $n-1$ to the system state $\vec{z}_n$ at time step $n$
\begin{equation}
    \vec{z}_n = \mathbf{U} \vec{z}_{n-1},
\end{equation}
where $\mathbf{U}$ is the monodromy matrix. The state vector $\vec{z}_n$ is a discrete representation of the DDE state, which contains information of the system variable $\vec{y}$ for the time interval $[t-\tau_{\text{max}},t]$, where $\tau_{\text{max}}$ is the maximum delay. The matrix $\mathbf{U}$ is a high dimensional approximation of the monodromy operator which is an operator that allows for mapping dynamic states forward in time by one period \cite{Khasawneh2013}. The full operator is infinite dimensional, but this method utilizes finite approximations of the operator to permit computing approximate solutions to the characteristic equation of the system. Specifically for the spectral element method, because the system is transformed into a discrete map, if the eigenvalues of $\mathbf{U}$ have a magnitude less than unity, that equilibrium point is stable and larger than one makes it unstable. If the magnitude is exactly one the equilibrium point is said to be marginally stable and this also is indicative that the equilibrium is non-hyperbolic \cite{Kuznetsov2004}. Note that the monodromy matrix for an autonomous system always has a trivial eigenvalue $\lambda=1$ which does not dictate the stability of the equilibrium point\cite{BEYN2002149}. In the case where the trivial eigenvalue is the furthest from the origin we take the second largest eigenvalue of the system to characterize the stability.

This problem falls under the broader classification of \textit{pseudospectral differencing methods} where in general an approximation of an infinite dimensional operator is computed resulting in a matrix where the eigenvalues approach solutions to the characteristic equation of the linear delay differential equations \cite{breda2005pseudospectral}. Further, Breda et al. proved many useful convergence results for such methods such as the fact that none of the eigenvalues of the approximation matrix are ``ghost roots". This means that all of the computed eigenvalues will eventually converge to a true root of the full infinite dimensional system if sufficient nodes are used in the discretization \cite{breda2005pseudospectral}. It is also known that roots closer to the origin in the complex plane are approximated first with these differencing methods so it is important to use sufficient discretization meshes to find the unstable eigenvalues \cite{breda2005pseudospectral}. It has been shown that for a DDE system the number of characteristic equation roots in the right half of the complex plane is finite \cite{stepan1989} meaning that if enough eigenvalues are approximated for the system about its equilibrium, eventually the right most eigenvalue will be computed which allows for characterizing the stability of the equilibrium point. For a discrete system this means that the number of eigenvalues outside of the unit circle is finite. In further sections, methods described in \cite{Khasawneh2013} are applied for discretizing and computing dominant eigenvalues for the metabolic systems to evaluate the system stability at different parameters. 

\subsection{Numerical Simulations}
 For the second method, we chose to perform many numerical simulations to demonstrate the behavior of the system and connect the results to experimental observations. This was done by brute force simulation of the system using the \texttt{Julia} differential equations library. The goal was to study specific features of the system trajectories in the parameter space to locate regions with different types of solutions. A diagram is obtained from the simulations by computing scalar features of the asymptotic solutions from time domain simulations and plotting the result as an image as a 2D projection of the overall parameter space. The features can be used to distinguish periodic solutions from equilibria. If the simulation times are long enough such that the system behaviour is characteristic of its long run behavior, we refer to this as the \textit{steady state response} as this is when the transient response has dissipated. The method for computing these response feature diagrams for this system was inspired by the \textit{AttractionsViaFeaturizing} function of the dynamical systems library in \texttt{Julia} \cite{Datseris2018}. This method computes a feature $M: \mathbb{R}^n\rightarrow \mathbb{R}$ on the system trajectory with $n$ system variables that indicate different features of the system response at each point in the parameter space. For this analysis, we needed to fix the history function for our systems. Specifically for this paper, we focus on a feature based on the amplitude of the time series signals. However, many other features can be used to study system behavior such as the mean response and standard deviation. We compute the amplitude feature $A$ of a response $x_i(t)$ as:
\begin{equation*}
    A_i = \frac{1}{2}\left(\max{(x_i(t)}-\min{(x_i(t))}\right).
\end{equation*}
The amplitude feature is then consolidated into a scalar value by summing over the variables in $i$ as:

\begin{equation}\label{eq:amp_feature}
    M_A=\sum_{i=1}^n A_i.
\end{equation}
If the trajectory is stable, we expect $M_A$ to be close to zero and if the response contains oscillations $M_A$ should be nonzero and finite.

\subsection{Low Growth History Functions}

To compute features of the system response, sufficient information is required to simulate the system such as the start time, end time and initial conditions. One critical difference between time delay differential equation systems (DDE) and ordinary differential equation systems (ODE) is that a DDE system requires the solution to be defined over the interval $[-\tau_{\text{max}},0]$ rather than just supplying a single point initial condition for an ODE system. A history function was chosen based on the experimental process for achieving these metabolic oscillations in practice \cite{Tu2005}. In this paper, the authors starve the cells of all resource prior to the oscillations. Consequently, the protein production rate is also zero during this time. 

\subsubsection{Single Protein History Function and Initial Conditions}\label{sec:1ptn_history}
 To define the history function for the single protein model, we assume that the cell was operating at zero protein production on $[-\tau,0)$. In other words, $p(\theta)=0$ and $R(\theta)=0$ for $\theta \in [-\tau,0)$, whereas at time $t=0$ we set $p(0)=p_0$. We obtain the initial conditions of the system by letting $t=0$ in Eq.~\eqref{eq:single_protein_resource} yielding:
\begin{equation}\label{eq:r0constraint}
	R_T=R_0+A\int_{-\tau}^0 f(p(s))R(s)ds,
\end{equation}
where $R_0=R(0)$ is the resource value at time $t=0$. For a response of this system to be valid, Eq.~\eqref{eq:r0constraint} must hold for the value of $R_T$ used for the simulation. Since $p$ and $R$ are zero for the history function, the integral in Eq.~\eqref{eq:r0constraint} vanishes and we have $R_0=R_T$. As a result, for each simulation $R_T$ can be specified and any positive value of $p_0$ can be chosen. This system can then be studied by varying the parameters $p_0$, $\tau$ and $R_T$ and holding the remaining parameters constant to determine which parameter values result in periodic solutions.

\subsubsection{Three Protein History Function and Initial Conditions}\label{sec:3ptn_history}
The solutions for the three protein system are also studied using the same zero resource and zero protein production assumption prior to $t=0$. So $p_1(\theta)=p_2(\theta)=p_3(\theta)=0$ and $R(\theta)=0$ for $\theta \in [-\tau_{\text{max}},0)$. Applying this to the total resource equation from Eq.~\eqref{eq:rt} again yields $R_0=R_T.$
Similar to the single protein system, the initial protein production rates $p_{10}$, $p_{20}$ and $p_{30}$ can be varied in the system. The dimension of the parameter space is reduced by limiting this system to the case where $\tau_1=\tau_2=\tau_3=\tau$ or in other words, all three proteins require the same production time. We then vary this delay and the total resource to determine which parameter combinations yield oscillations in the system response.

\subsection{Boundary Value Calculation of Periodic Solutions to Nonlinear DDE Systems} \label{sec:bvp_sol}
As an alternative to detecting periodic orbits by simulation and amplitude computation, we will find them directly by  solving a  boundary value problem (BVP). This method comes from \cite{khasawneh2012periodic} where the authors describe how a nonlinear DDE system can be converted to a BVP where the solutions to this problem correspond to periodic solutions of the original system. Specifically, the system is converted to the boundary value problem in Eq.~\eqref{eq:BVP} and the spectral element method is used to discretize the DDE system to approximate the infinite dimensional BVP as a finite dimensional problem that can be solved numerically to obtain a single period of the periodic solution to the system \cite{khasawneh2012periodic}. 

\begin{equation} \label{eq:BVP}
\begin{split}
    \vec{f}&=\frac{d\vec{x}}{dt}-T\vec{g}(\vec{x}(t),\vec{x}(t-\tau/T))=0,~ t\in[0,1],\\
    \vec{x}&(s) - \vec{x}(s+1) = 0,~
    s\in[-\tau/T,0],\\
    p&(\vec{x})=0,
\end{split}
\end{equation}
where $T$ is the system period, $\vec{x}$ is the vector of system variables, and the first line corresponds to the specific DDE system being studied. In this case we take $\vec{g}$ to be Eq.~\eqref{eq:g_single} for the single protein system, and Eq.~\eqref{eq:g_triple} for the three protein system. The second line imposes a periodicity condition on the system and the last line imposes a phase condition to yield a unique periodic solution by setting $p(\vec{x})$ to be the inner product of the initial state ($\vec{x}_0$) and the time derivative $\dot{\vec{x}}(t)$ \cite{khasawneh2012periodic}. An initial guess is provided by simulating the system in \texttt{Julia} using the differential equations library and this simulation is provided to the boundary value problem solver in \texttt{Matlab} to perform Newton-Raphson iteration and converge on the periodic solution. This process also yields an approximation to the period $T$ of the system. This method has been shown to compute accurate periodic solutions for nonlinear DDE systems with exponential convergence rates as the number of mesh points increases \cite{khasawneh2012periodic} making it ideal for verifying parameters of the metabolic system that result in oscillating solutions to the system. 

%% file: Sections/results.tex
\section{Results - Single Protein}\label{sec:1ptnresults}

The single protein system is analyzed in this section by applying the three methods outlined in Sec.~\ref{sec:methods}.

\subsection{Spectral Element Linear Stability}
The eigenvalues of the linearized single protein system about its nontrivial equilibrium points were approximated at each set of parameters in a $400\times400$ grid in the $(\tau, R_T)$ plane varying each parameter from zero to 50 using the monodromy matrix from the spectral element method \cite{Khasawneh2013}. We hold the remaining parameters constant at $\kappa=0.5$, $A=1.0$, $B=2.0$, $D=10.0$, $n=2$. The monodromy matrix requires an oscillation period to map the system states to the next period. It has been shown that for systems with one delay the period can be set as the delay for stability computations \cite{Khasawneh2013}. For this reason, we set the period to $\tau$ for this stability analysis. We examine the stability of the equilibrium points by plotting the magnitude of the eigenvalue furthest from the origin. Stability diagrams were plotted for the nontrivial equilibrium points. Note that below the curve in Eq.~\eqref{eq:saddlenode}, only the trivial equilibrium is present $(0,R_T)$, but above this curve three equilibrium points exist in the system. We only consider the stability of the nontrivial equilibria in this section as the stability of the trivial solution is computed analytically in Section~\ref{sec:single_lin}. As a result, we color points in the stability diagram as white if only the trivial equilibrium is present in that region.

\begin{figure}[t]
    \centering
    \includegraphics[width=0.9\textwidth]{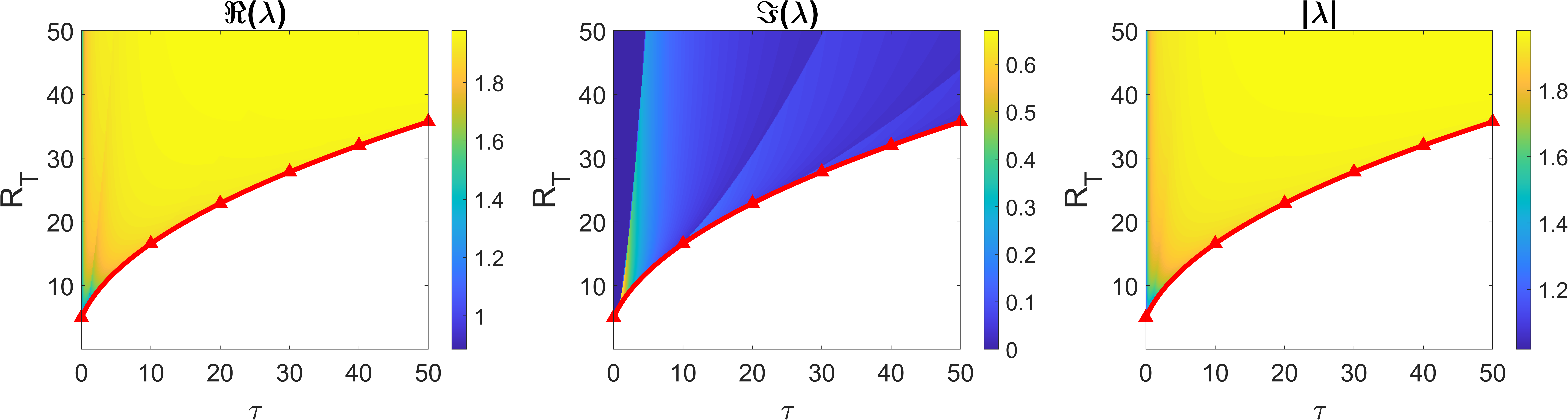}
    \caption{Single protein middle equilibrium point stability diagrams. Specifically, the eigenvalues with maximum magnitude of the monodromy matrix are plotted with respect to the parameters $\tau$ and $R_T$. (left) the real part of the dominant eigenvalue, (middle) imaginary part of the dominant eigenvalue, (right) the modulus of the eigenvalue. The red curve with triangles is the saddle node boundary that separates regions with 1 and 3 equilibria. Above the red curve all three equilibrium points exist and below the curve only the trivial point is present.}
    \label{fig:middle_stability}
\end{figure}

We start by computing the stability of the middle equilibrium point as shown in Fig.~\ref{fig:middle_stability} where we plot the dominant eigenvalue of the middle equilibrium point for combinations of $\tau$ and $R_T$ between 0 and 50. We see that for all parameters shown, this equilibrium point is unstable because its largest eigenvalue is outside of the unit circle in the complex plane. 

\begin{figure}[htbp]
    \centering
    \includegraphics[width=0.9\textwidth]{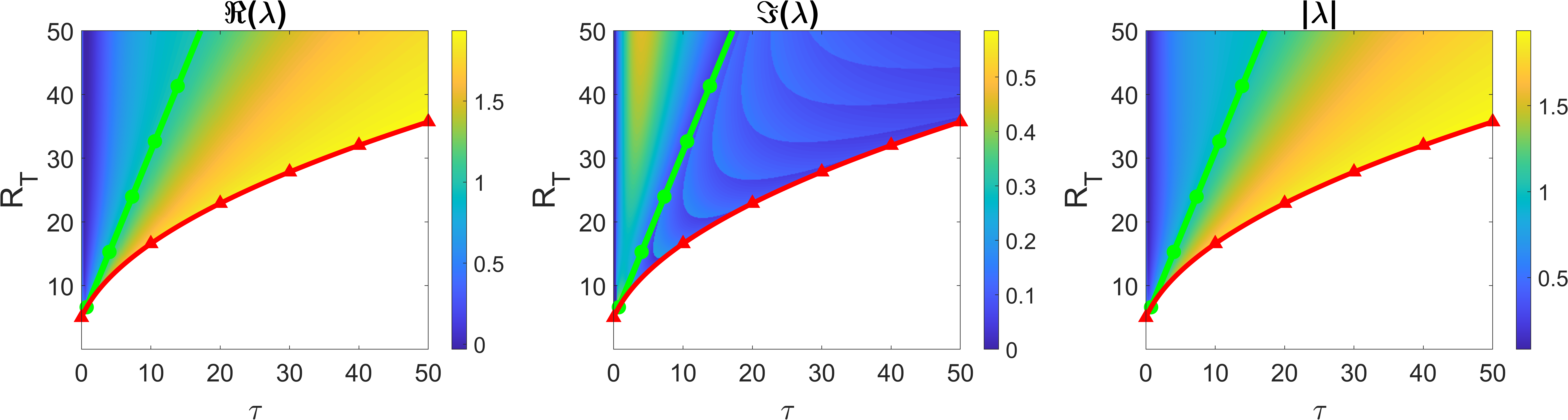}
    \caption{Single protein top equilibrium point stability diagrams. Specifically, the eigenvalues with maximum magnitude of the monodromy matrix are plotted with respect to the parameters $\tau$ and $R_T$. (left) the real part of the dominant eigenvalue, (middle) imaginary part of the dominant eigenvalue, (right) the modulus of the eigenvalue. The red curve with triangles is the saddle node boundary that separates regions with 1 and 3 equilibria. Above the red curve all three equilibrium points exist and below the curve only the trivial point is present. The Hopf bifurcation curve is shown as a line with green dots.}
    \label{fig:top_stability}
\end{figure}

Next, we plot the largest magnitude eigenvalue of the top equilibrium point in Fig.~\ref{fig:top_stability} where we see that for small delay and sufficient resource, the top equilibrium point is  stable with $|\lambda| < 1$.  As the delay increases for a given total resource, a pair of complex conjugate eigenvalues leave the unit circle that govern the stability of this equilibrium point making it an unstable focus \cite{Kuznetsov2004}.

Therefore, a Hopf bifurcation occurs from the top equilibrium point along this line. We plot the Hopf bifurcation curve in subsequent stability diagrams as a green line with dots. This line was found to be approximately,
\begin{equation}\label{eq:1ptn_hopf}
    R_T = 2.6449 \tau + 4.6323,
\end{equation}
for $\tau \geq 0.75$ by computing a linear regression along the boundary where the eigenvalue exits the unit circle. The model had a coefficient of determination of 0.9999 indicating that this boundary is well approximated by a linear model.

\subsection{Response Features}

\begin{figure}[htbp]
	\centering
	\begin{minipage}{0.45\textwidth}
		\centering
		\includegraphics[width=\textwidth]{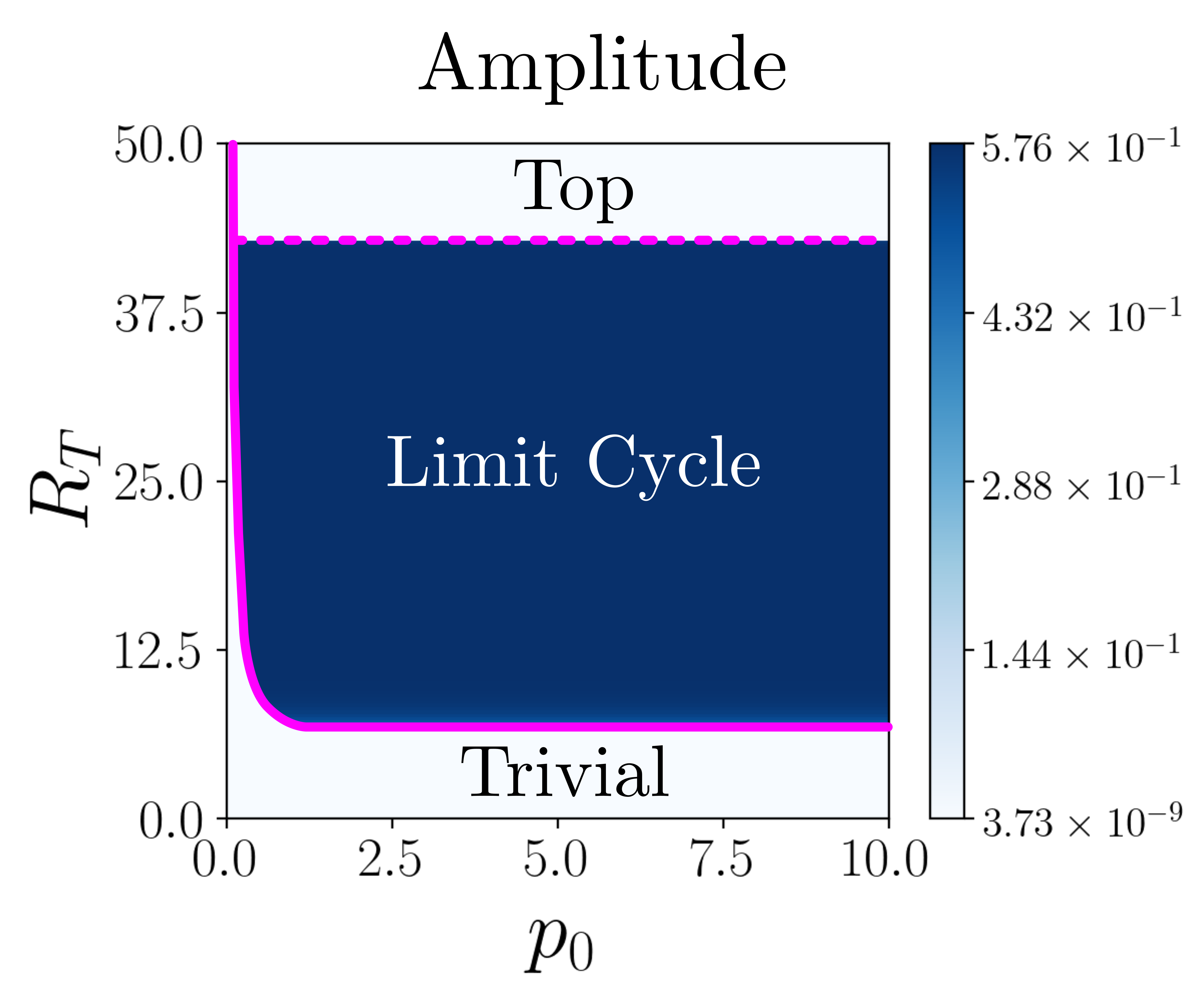}
		(a) $p_0$ Dependence
	\end{minipage}
	\begin{minipage}{0.45\textwidth}
		\centering
		\includegraphics[width=\textwidth]{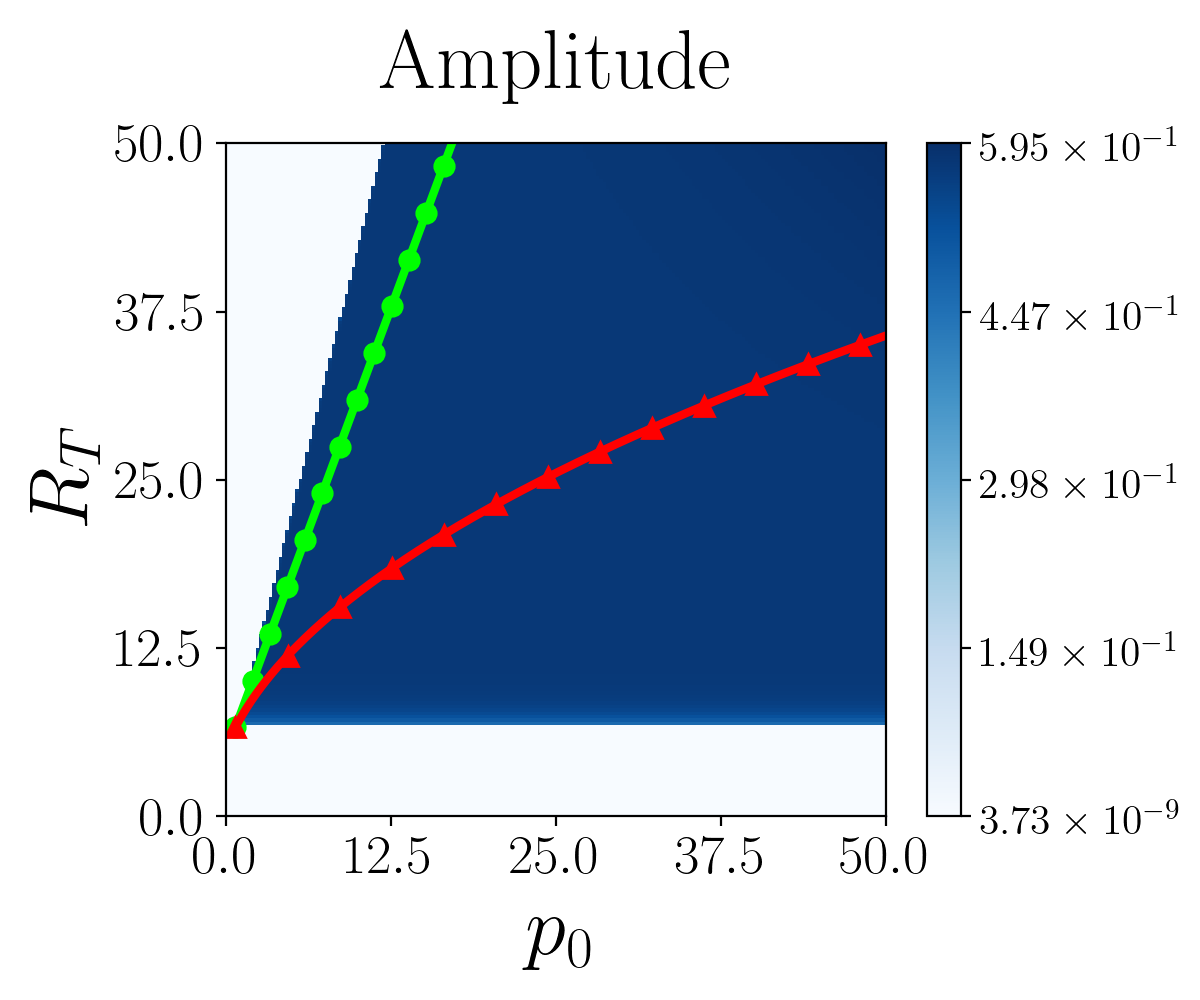}
		(b) $\tau-R_T$ Dependence
	\end{minipage}
	\caption{Single protein model response feature diagrams by varying system parameters $p_0$, $\tau$ and $R_T$ and simulating the system at for each parameter combination between 0 and 50. The solid pink curve corresponds to a point on the horizontal boundary at $\tau=25$ in the right image and the dashed pink curve corresponds to a point on the boundary above the green curve at $\tau=25$ in the right image. The red curve with triangles is the saddle node boundary that separates regions with 1 and 3 equilibria, the line with green dots is the Hopf bifurcation boundary.}
	\label{fig:1ptn_features}
\end{figure}

Response feature diagrams were generated for the single protein system (Eq.~\eqref{eq:single_protein_system}) for $\kappa=0.5$, $A=1.0$, $B=2.0$, $D=10.0$, $n=2$ with varying $\tau$, $R_T$ and $p_0$. The results for these simulations are shown in Fig.~\ref{fig:1ptn_features} where we color pixels in the parameter space according to the response amplitude feature using Eq.~\ref{eq:amp_feature}. We used the starving cell history function from Section~\ref{sec:1ptn_history} and each simulation was taken between 10000-11000 time units to ensure that the transient response had dissipated. The authors acknowledge the arbitrarily chosen parameters for this system and that these parameters may not be in biologically significant range. However, our model is conceptual and the purpose of this paper is to demonstrate that certain parameters yield oscillations in the protein production when the cell is starved of resource prior to $t=0$. This is also the reason why we use ``time units" instead of seconds for the simulations. 

First we study the dependence on the initial protein production rate $p_0$ by fixing the delay at $\tau=10$ and plotting the amplitude feature over the region $(p_0, R_T)\in[0,10]\times[0,50]$. We see in Fig.~\ref{fig:1ptn_features}(a) that for nontrivial $p_0$, the response is essentially independent of the initial condition so any large enough initial protein production rate was sufficient. For small $p_0$ the response approaches the trivial equilibrium point. While this diagram is only shown for a single delay, we observed a trend where as the delay varies, the only change is in the width of the limit cycle region for large enough $p_0$. For this reason, we arbitrarily choose $p_0=10$ for our initial $p_0$.

 Next, we keep $p_0=10$ and vary the parameters $(\tau, R_T)\in[0,50]\times[0,50]$ and plot the amplitude feature in this region of the parameter space in Fig.~\ref{fig:1ptn_features}(b) along with the Hopf and saddle node bifurcation boundaries obtained from the linear stability analysis. We see that periodic solutions were found above the Hopf curve for this particular history function indicating that the Hopf bifurcation is subcritical. So slightly above the green curve we have a bistability between the top equilibrium, trivial equilibrium, and the limit cycle. Below the Hopf curve we have a bistability between the trivial equilibrium point and the limit cycle and below $R_T\approx7$ we did not observe any oscillations and the trajectory approached the trivial equilibrium. Note that the pink curves in Fig.~\ref{fig:1ptn_features}(a) are specific to $\tau=10$ and will increase as the delay is increased according to the bounds of the periodic region in Fig.~\ref{fig:1ptn_features}(b). In other words, at a delay of 10 if we draw a vertical line in Fig.~\ref{fig:1ptn_features}(b), we should expect it to intersect the blue region at $R_T\approx7$ and $R_T\approx42$ which correspond to the pink curves in Fig.~\ref{fig:1ptn_features}(a) for nontrivial $p_0$.

\subsection{Periodic Solutions}
\begin{figure}[t]
	\centering
	\includegraphics[width=0.5\textwidth]{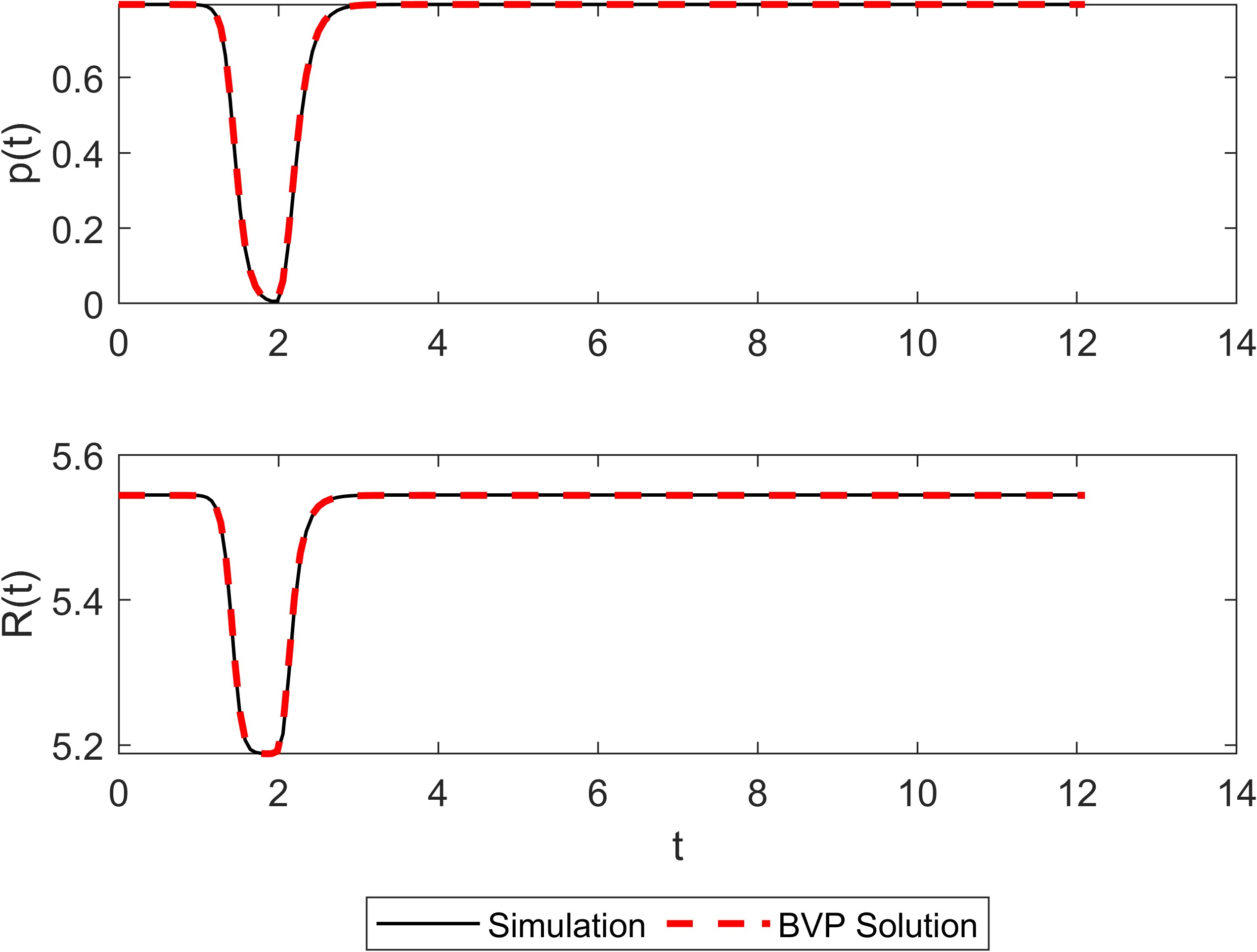}
	\caption{$\tau=12$ and $R_T=50$ single protein periodic solution results from solving the relevant boundary value problem. }
	\label{fig:1250periodic_sol}
\end{figure}

Next we utilize the spectral element approach to solve the boundary value problem in Eq.~\eqref{eq:BVP}. This process was performed on the three points in the $(\tau, R_T)$ parameter space where oscillations were expected and the two points where we expect fixed point responses. The first point considered was $\tau=12$ and $R_T=50$. We see that this point corresponds to a response with nonzero amplitude indicating that oscillations should be expected and is in the subcritical region of the Hopf bifurcation. The system was simulated and sampled between 15,988-16,000 time units for the period. Because the system is autonomous, we can take the period to be equal to the delay.  Solving the boundary value problem in Eq.\eqref{eq:BVP} for the periodic solution, we obtain the response shown in Fig.~\ref{fig:1250periodic_sol}. We see that the periodic solution from the boundary value problem closely matches the simulation result with the period matching the delay. Further, the protein production rate is nearly constant and close to the top equilibrium point $(p^*, R^*)=(0.7434, 5.3982)$ for most of the period in this case with a drop in the production rate emerging yielding the metabolic oscillations.

\begin{figure}[htbp]
	\centering
	\includegraphics[width=0.5\textwidth]{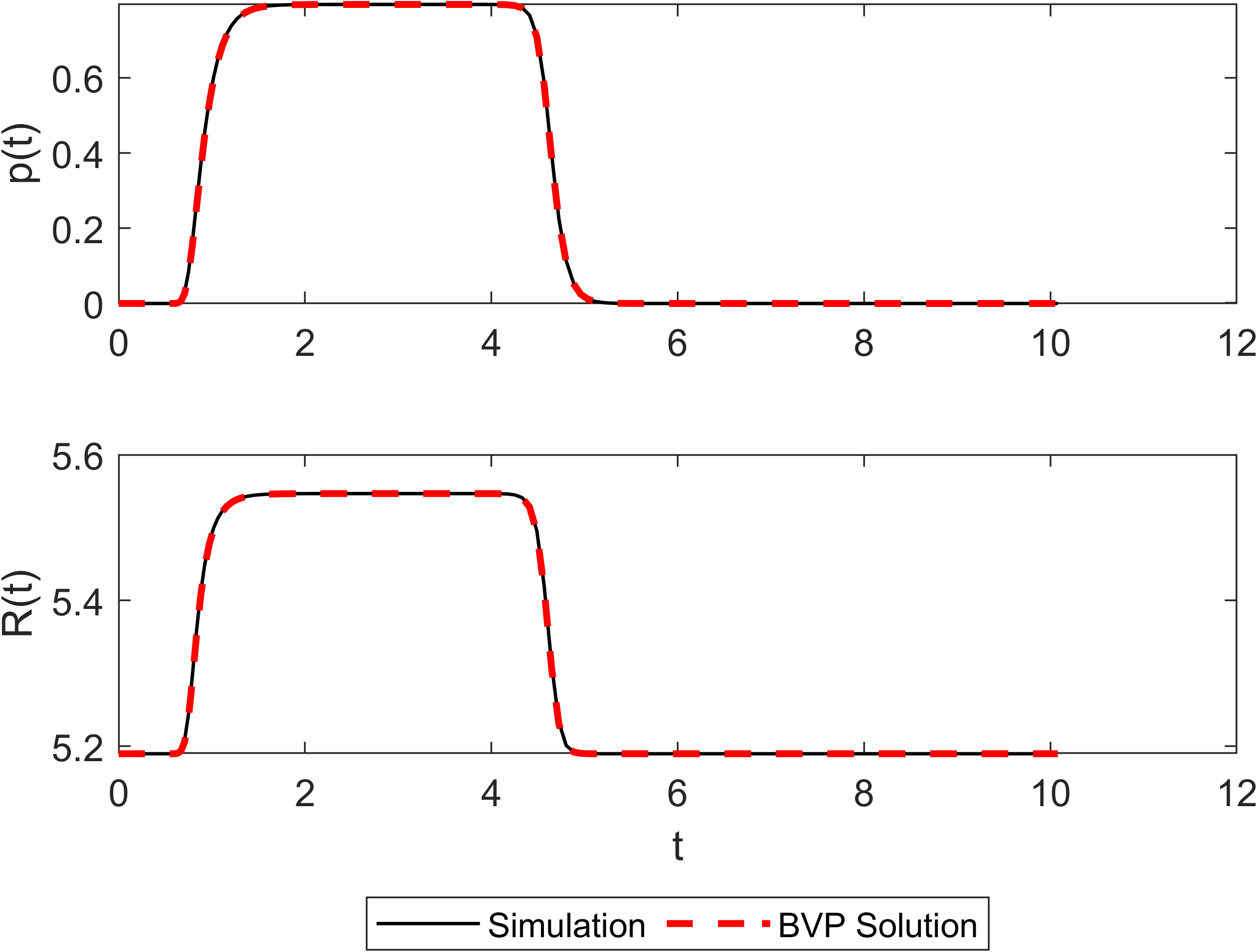}
	\caption{$\tau=10$ and $R_T=20$ single protein periodic solution results from solving the relevant boundary value problem. }
	\label{fig:1020periodic_sol}
\end{figure}

The next parameters that were considered were $\tau=10$ and $R_T=20$. The system was simulated at these parameters from 15,990-16,000 time units and the period was set to 10. Passing this initial guess into the boundary value problem, the obtained periodic solution is shown in Fig.~\ref{fig:1020periodic_sol}. Interestingly, we see that the time that the protein production rate spends at 0 is much longer compared to Fig.~\ref{fig:1250periodic_sol}. As the Hopf bifurcation curve is crossed, the drop in the protein production rate appears to spend more time at zero during the oscillation period. 

\begin{figure}[htbp]
	\centering
	\includegraphics[width=0.5\textwidth]{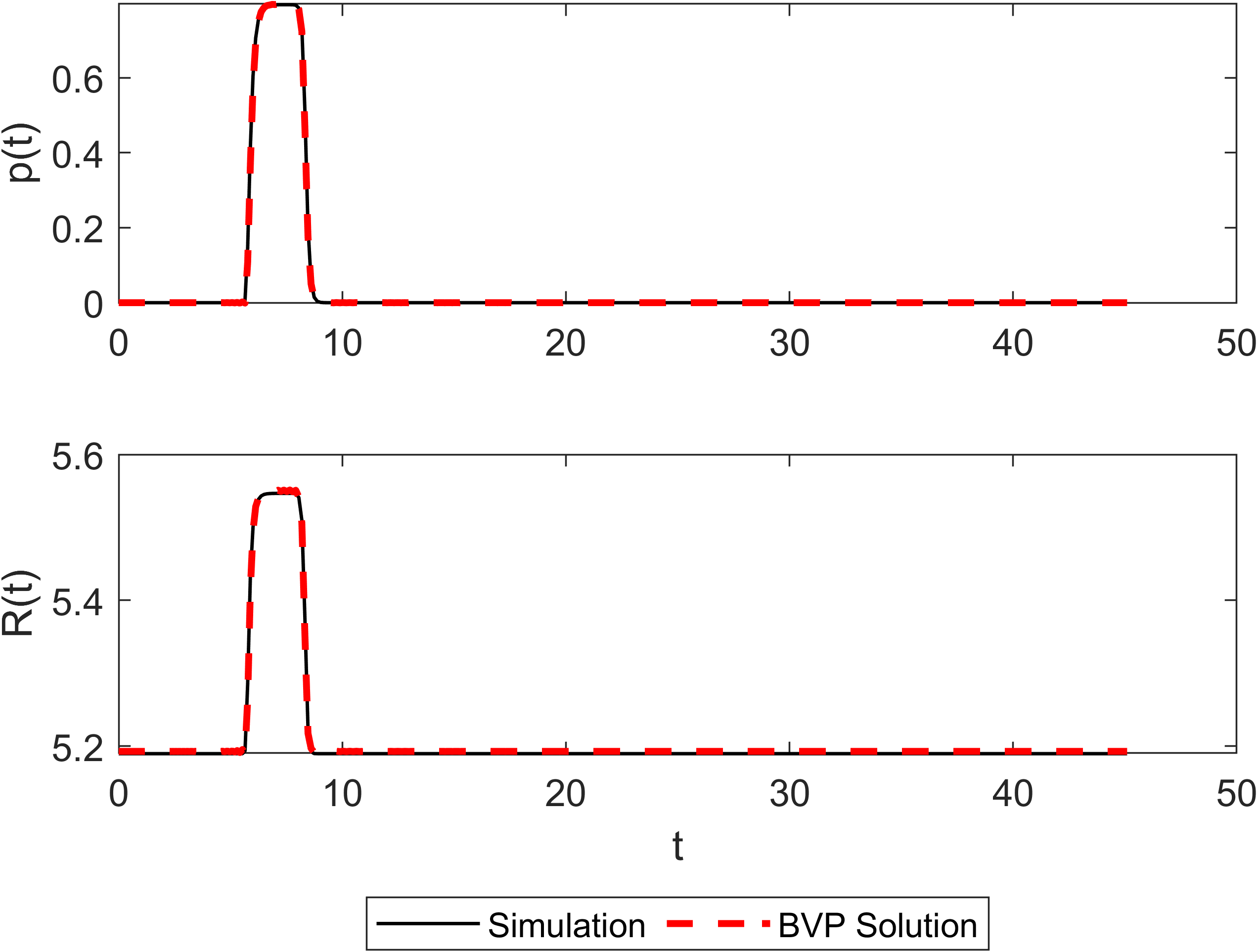}
	\caption{$\tau=45$ and $R_T=15$ single protein periodic solution results from solving the relevant boundary value problem. }
	\label{fig:4515periodic_sol}
\end{figure}

The third set of parameters considered was $\tau=45$ and $R_T=15$. The system was simulated at these parameters from 15,955-16,000 time units and the period was set to 45 resulting in the periodic solution shown in Fig.~\ref{fig:4515periodic_sol}. We see that the trajectory starts to spend more time near a protein production rate of zero as the total resource approaches the horizontal line $R_T\approx7$ in Fig.~\ref{fig:1ptn_features}(b). The periodic solutions shown demonstrate the transition from the fixed point stability at the top equilibrium to fixed point stability at the trivial equilibrium.

\subsection{Steady State Solutions}
\begin{figure}[htbp]
	\centering
	\includegraphics[width=0.6\textwidth]{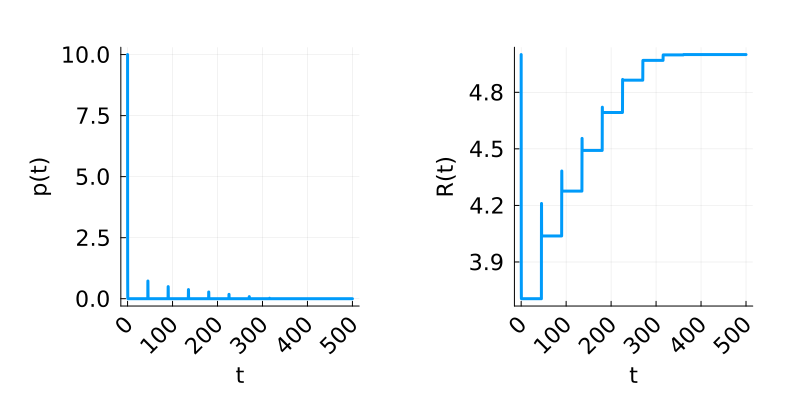}
	\caption{Approach to trivial fixed point  for low growth conditions in the cell ($\tau=45$, $R_T=5$)}
	\label{fig:4505resp}
\end{figure}

\begin{figure}[htbp]
	\centering
	\includegraphics[width=0.6\textwidth]{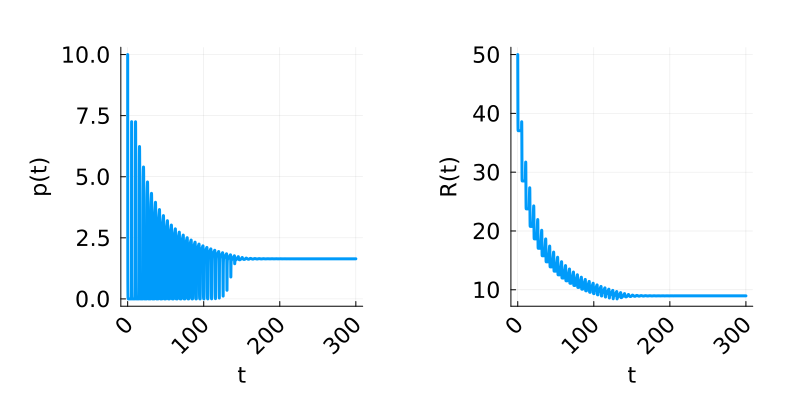}
	\caption{Fixed point response for high growth conditions in the cell ($\tau=5$, $R_T=50$)}
	\label{fig:0550resp}
\end{figure}

We also explore the steady state solutions of the system by examining two parameter conditions. The first case is where the total resource is too low to sustain protein production (low growth conditions). In this case, we found a trajectory that approaches the trivial equilibrium point. This was verified by simulating the system at $\tau=45$ and $R_T=5$. The resulting response is shown in Fig.~\ref{fig:4505resp} where we see the system approach $(p,R)=(0,R_T)$.
We also consider the case where the cell has access to plentiful resources and can synthesize proteins at a constant rate (high growth conditions). To examine this case, we simulated the system at $\tau=5$ and $R_T=50$. The response for these parameters is shown in Fig.~\ref{fig:0550resp}. We see that as time progresses, the protein production rate approaches a steady state value because the cell is able to produce proteins at a constant rate. Further, the point that this trajectory approaches corresponds to the top equilibrium point of the system which for these parameters works out to be $(p^*,R^*)\approx(1.6412, 8.968)$.

\subsection{Single Protein Summary}

Three distinct behaviors were observed in the single protein time delay model. First, if the resources are not sufficient to sustain metabolic activity, the system will approach the trivial equilibrium point with zero protein production rate. If the resources are plentiful, they can sustain constant protein production at the top equilibrium point. Between these two cases, the cell initially has enough resource to synthesize proteins, but as the protein production rate increases, resources are used and the metabolic activity decreases. This balance between constant production and no production seems to lead to oscillations in the system response. Slightly above the Hopf bifurcation curve, we observe oscillations that are close to a constant solution at the top equilibrium and as the parameters cross the Hopf curve and approach the line $R_T\approx7$, the solution continues to oscillate but with a solution that is closer to a constant solution at the trivial equilibrium. The oscillation region represents a transition between the top and trivial equilibria and the middle solution remains unstable for all parameters in this region. For $\tau<0.75$, the solution can switch directly between the top and trivial equilibrium with no oscillations, but after this bifurcation point at $(\tau, R_T)\approx(0.75, 6.6)$, the periodic solution emerges to transition from constant to zero protein production.

\section{Results - Three Proteins}\label{sec:3ptnresults}
The three protein system is analyzed in this section by applying the three methods outlined in Sec.~\ref{sec:methods}.

\subsection{Spectral Element Linear Stability}\label{sec:3ptn_lin_stability}
Our main goal with the three protein system was to find parameters where the protein production rates peak at different times in the period. This phenomena would be indicative of the cell prioritizing its resources to produce proteins in a way that could be more efficient. To begin exploring this systems parameter space, we use the spectral element method to study the linear stability of the three protein system with arbitrarily chosen parameters $\kappa=0.5$, $A=1.0$, $B_1=2.0$, $B_2=2.0$, $B_3=2.0$, $D_1=10.0$, $D_2=10.0$, $D_3=10.0$, $n=2$. We set $\tau_1=\tau_2=\tau_3=\tau$ such that all three proteins require the same amount of production time, and vary $\tau$ and $R_T$ just as was done with the single protein system.

However, for this system we do not have analytical expressions for the equilibrium solutions and we only have the coupled polynomial system in Eqs.~\eqref{eq:poly_eq_1}, \eqref{eq:poly_eq_2}, and \eqref{eq:poly_eq_3}. Solving this system of equations is a nontrivial task, but if we make some assumptions based on our observations from the single protein system we can still generate stability diagrams using this method. Namely, we will assume that this system also exhibits three possible equilibrium points (top, middle, and trivial). Using the variable precision accuracy (VPA) solver in \texttt{Matlab}, we can solve this system of equations numerically in our parameter space and approximate the dominant eigenvalues to characterize the stability of each point. The documentation for the VPA solver used states that for polynomial systems, all solutions in a region will be returned by the function \cite{vpasolve}. This solver was applied to a 400$\times$400 grid of parameters in $(\tau, R_T)\in[0,50]\times[0,100]$, and the assumption was found to be correct where one region of the space contained 3 equilibrium points and the other region only contained the trivial point.  

With the single protein system, we defined the top and middle equilibrium points based on the magnitude of the equilibrium protein production rate $p^*$. However, in the three protein system we have multiple equilibrium protein production rates. To modify this approach for the three protein system, we form the following vector of equilibrium coordinates,
\begin{equation*}
    \vec{p}=\begin{bmatrix}
        p_1^* & p_2^* & p_3^*
    \end{bmatrix}^T.
\end{equation*}

Thus, the top and middle equilibria are defined by the $l_2$ norm of $\vec{p}$ where the top equilibrium has the largest $l_2$ norm and the middle solution has a norm between the top and trivial. We then use Eq.~\eqref{eq:r_star} to obtain $R^*$ for a given set of parameters at each equilibrium point. This system can then be linearized about each equilibrium point using Eq.~\eqref{eq:3ptn_linear} and the dominant eigenvalue of the linearized system at a given set of parameters can be approximated with the spectral element method described in Section~\ref{sec:spec_stability}. We note that because a single delay is present, we use this delay for the system period when computing the monodromy matrix for the system.

\begin{figure}[t]
    \centering
    \includegraphics[width=0.9\textwidth]{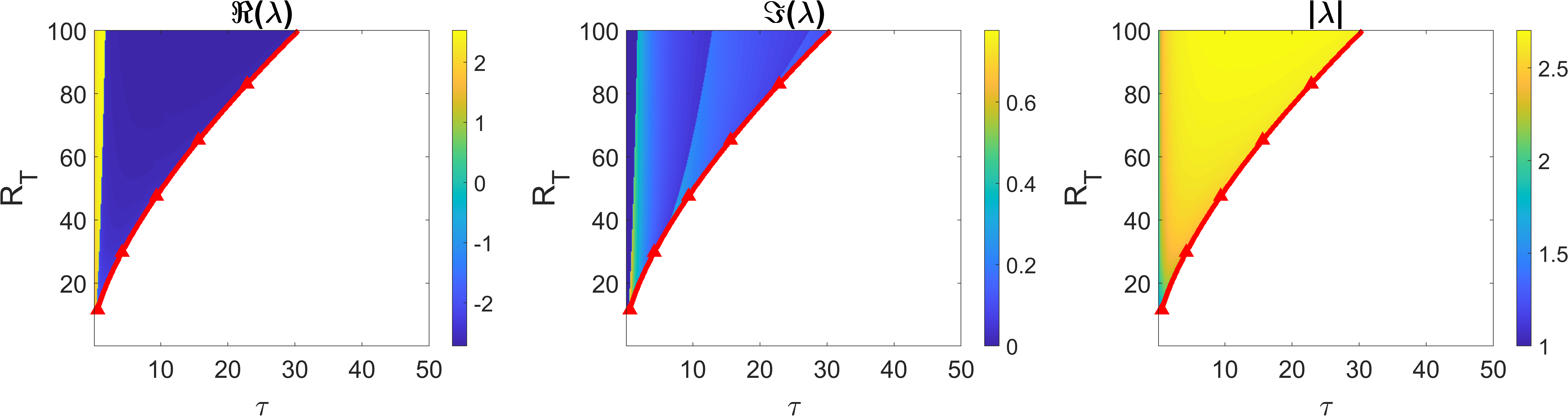}
    \caption{Three protein middle equilibrium point stability diagrams at equal delay. Specifically, the eigenvalues with maximum magnitude of the monodromy matrix are plotted with respect to the parameters $\tau$ and $R_T$. (left) the real part of the dominant eigenvalue, (middle) imaginary part of the dominant eigenvalue, (right) the modulus of the eigenvalue. The red curve with triangles is the saddle node boundary that separates regions with 1 and 3 equilibria. Above the red curve all three equilibrium points exist and below the curve only the trivial point is present.}
    \label{fig:3mid_evs}
\end{figure}

Similar to the single protein model, we start by plotting the stability of the middle equilibrium point in Fig.~\ref{fig:3mid_evs}. Note that we color the region as white if only the trivial equilibrium is present. We see that there are two distinct regions in the plot of the magnitude of the eigenvalue. The curve that separates these regions is the saddle node bifurcation curve for this system. Using a third order polynomial fit, the saddle node curve was found to be,

\begin{equation}
    R_T=0.0016\tau^3 -0.1118\tau^2 + 4.8855\tau +10.3749,
\end{equation}
for $\tau\geq0.625$. This curve had a coefficient of determination of 0.9997. We plot the saddle node boundary as a red curve with triangles in the stability diagrams. 

Because the dominant eigenvalue for the middle equilibrium always has a modulus greater than one for these parameters, this equilibrium point will not govern the stability of the system if another point is stable or marginally stable. This was also the case with the single protein system.

\begin{figure}[htbp]
    \centering
    \includegraphics[width=0.9\textwidth]{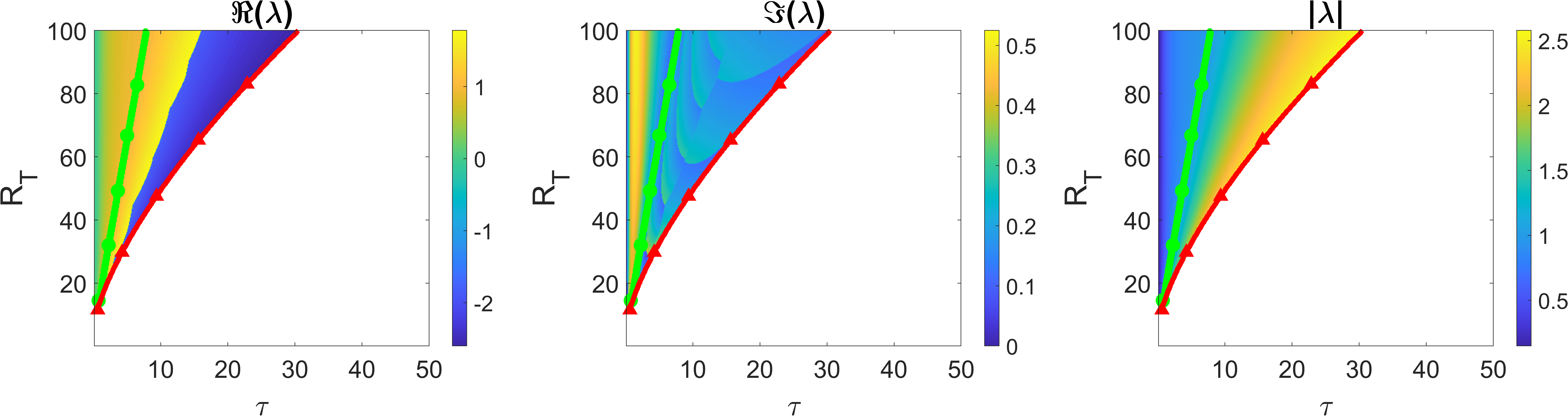}
    \caption{Three protein top equilibrium point stability diagrams at equal delay. Specifically, the eigenvalues with maximum magnitude of the monodromy matrix are plotted with respect to the parameters $\tau$ and $R_T$. (left) the real part of the dominant eigenvalue, (middle) imaginary part of the dominant eigenvalue, (right) the modulus of the eigenvalue. The red curve with triangles is the saddle node boundary that separates regions with 1 and 3 equilibria. Above the red curve all three equilibrium points exist and below the curve only the trivial point is present. The Hopf bifurcation curve is shown as a line with green dots.}
    \label{fig:3top_evs}
\end{figure}

Lastly, we plot the stability of the top equilibrium point of this system in Fig.~\ref{fig:3top_evs}. These diagrams show that for small delay and large resource, this point is stable. As the delay increases for a given resource, this point becomes an unstable focus by way of a Hopf bifurcation. We plot the Hopf bifurcation curve as a green line with dots. This curve was approximated by locating points in the parameter space with unit length eigenvalues. Using linear regression, the Hopf bifurcation curve was approximated to be,
\begin{equation}
    R_T=12.0948\tau + 4.7910,
\end{equation}
for $\tau \geq 0.75$. This model had a coefficient of determination of 0.9987 suggesting that it is a good approximation. Note that while there is interesting behavior that occurs between the green (dots) and red (triangles) curves in these stability diagrams, all of the eigenvalues plotted are outside of the unit circle and are therefore unstable so this is not a bifurcation it just means that another eigenvalue moved further from the origin. 

\subsection{Response Features}
\begin{figure}[t]
    \centering
    \includegraphics[width=0.5\textwidth]{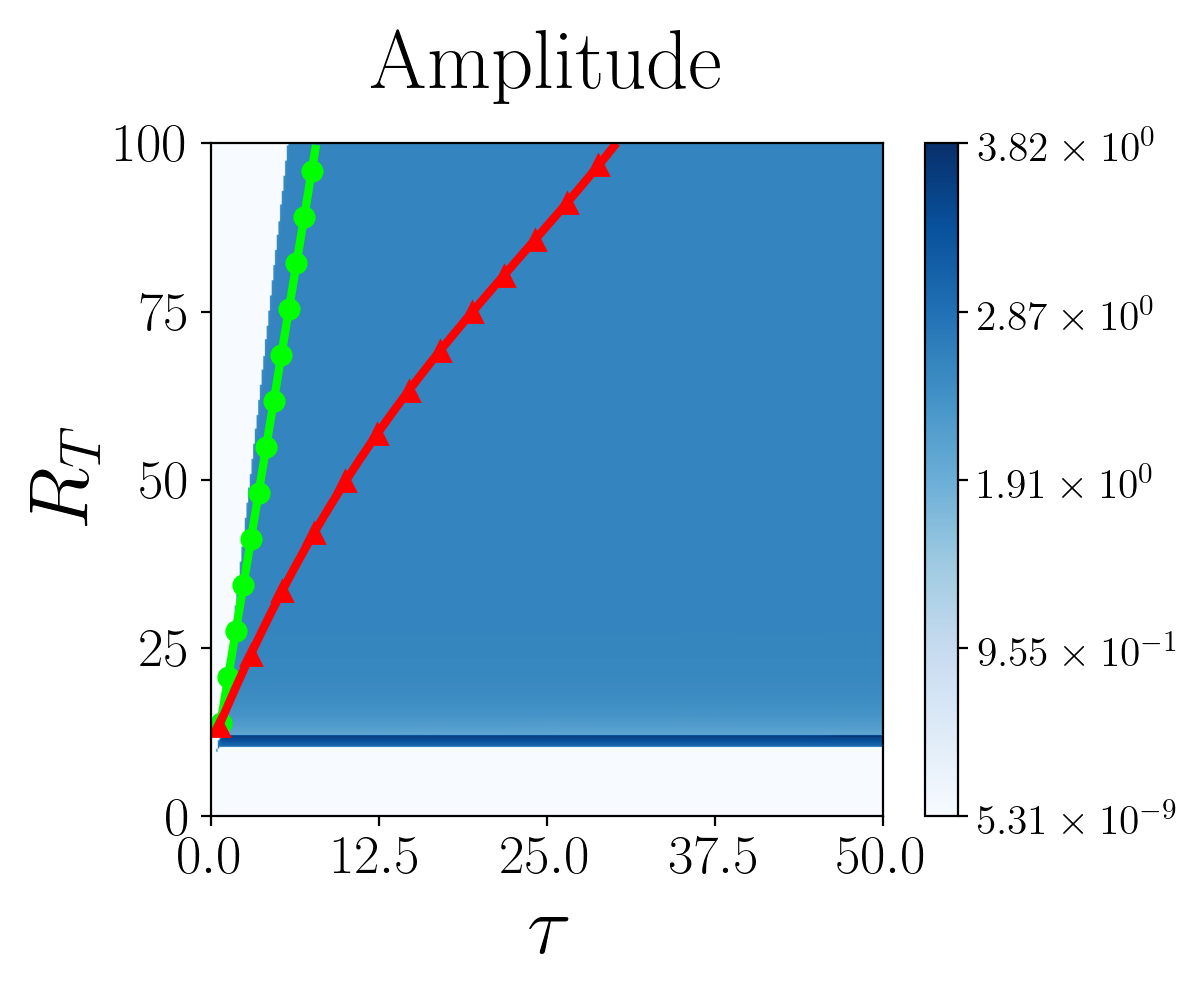}
    \caption{Three protein system response amplitude diagram in the $\tau-R_T$ parameter space where $\tau$ is the same delay for all three proteins. The low growth history function was used for all simulations in this diagram.}
    \label{fig:3ptn_amplitude}
\end{figure}

Holding the remaining parameters constant at the values from Section~\ref{sec:3ptn_lin_stability}, an amplitude feature diagram was generated with equal delays for all three proteins and varying the delay $\tau$ with the total resource $R_T$. The three protein system was simulated between 10,000-11,000 time units using the zero history function described in Sec.~\ref{sec:3ptn_history}, and the amplitude feature was plotted in the $\tau-R_T$ space where $\tau$ is the same for all three proteins. The results are shown in Fig.~\ref{fig:3ptn_amplitude}. We see that the amplitude diagram has a similar structure to the single protein system, but there is a change in the amplitude feature for a small region at low total resource. System responses in these regions are explored further in Section~\ref{sec:psol_eqdelay}. We plot the Hopf bifurcation curve from Section~\ref{sec:3ptn_lin_stability} in Fig.~\ref{fig:3ptn_amplitude} as the green line with dots. The red curve with triangles corresponds to the saddle node bifurcation curve from the approximations in  Section~\ref{sec:3ptn_lin_stability}. 

\subsection{Periodic Solutions}\label{sec:psol_eqdelay}

 Three points were considered within the region of the parameter space with nontrivial amplitude in Fig.~\ref{fig:3ptn_amplitude}. Namely, we choose $(\tau, R_T)=(5.7,100)$, $(25, 50)$ and $(25, 11.8)$. The first point is in the nonzero amplitude region to the left of the Hopf bifurcation curve. This point was chosen to verify the subcriticality of the Hopf bifurcation. The second point was chosen to show the response near the middle of the oscillation region. We chose the third point in the region of differing amplitude at low total resource to observe the changes that occur when the cell has limited resources available. The steady state regions or regions with near zero amplitude were found to exhibit similar behaviors to the single protein model outside of the oscillation region so these responses will not be considered. 

\begin{figure}[htbp]
    \centering
    \includegraphics[width=0.5\textwidth]{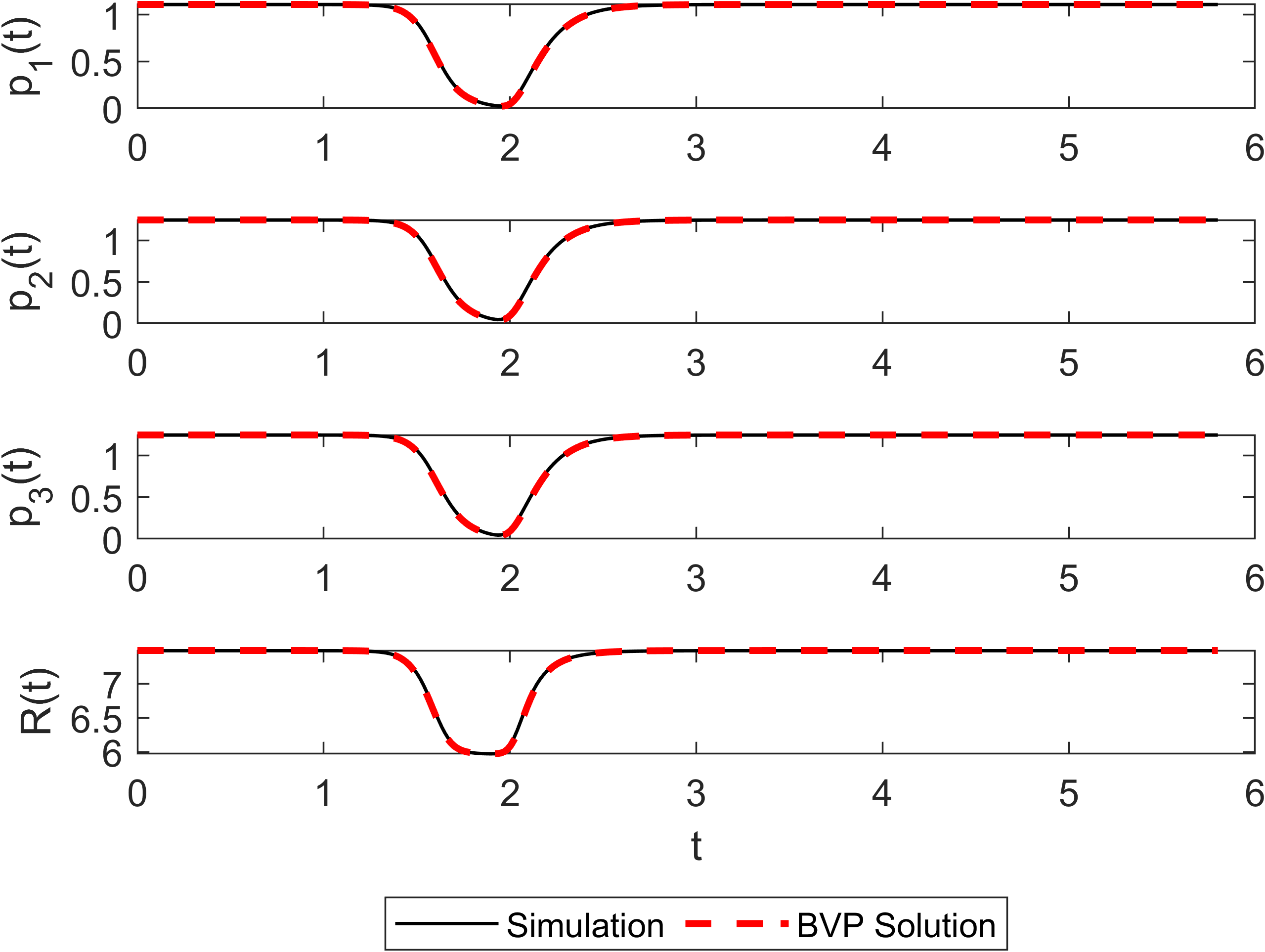}
    \caption{Boundary value problem solution for the three protein system with $\tau=5.7$ and $R_T=100$. }
    \label{fig:sametau_period57}
\end{figure}

Starting with $\tau=5.7$ and $R_T=100$, the system was simulated between 15,000 and 15,005.7 time units to capture a single period of the response. This solution was then verified by solving the nonlinear DDE boundary value problem to obtain the periodic solution in Fig.~\ref{fig:sametau_period57}. We see that the obtained solution is nearly identical to the simulation and appears to be close to a constant solution at the top equilibrium point which for these parameters is at $(p_1^*,p_2^*,p_3^*,R^*)\approx(0.9942,1.1328, 1.1328, 7.0965).$ We also note that all of the protein production rates here appear to oscillate in-phase.

\begin{figure}[htbp]
    \centering
    \includegraphics[width=0.5\textwidth]{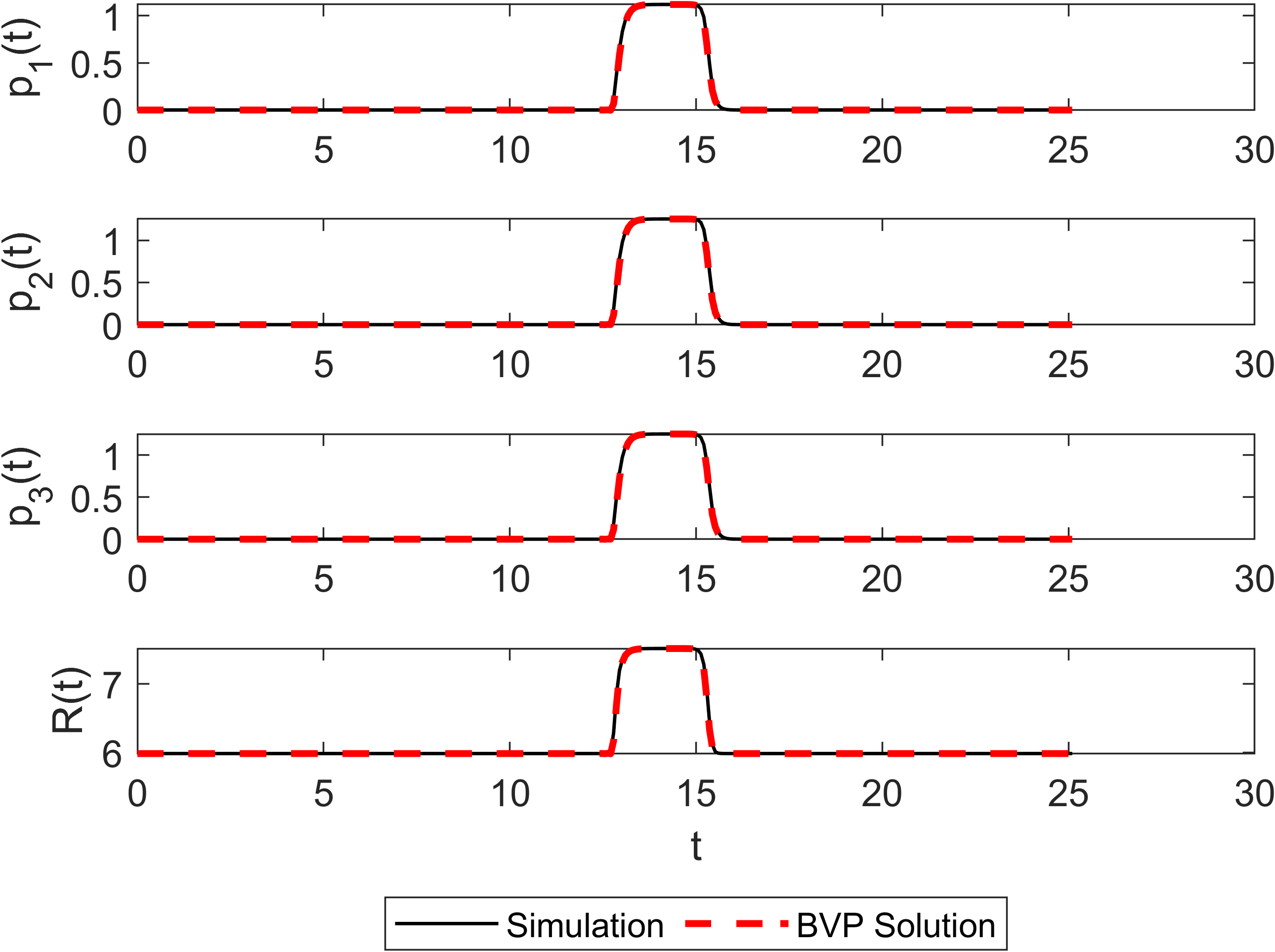}
    \caption{Boundary value problem solution for the three protein system with $\tau=25$ and $R_T=50$. }
    \label{fig:sametau_period2550}
\end{figure}

\begin{figure}[htbp]
    \centering
    \includegraphics[width=0.5\textwidth]{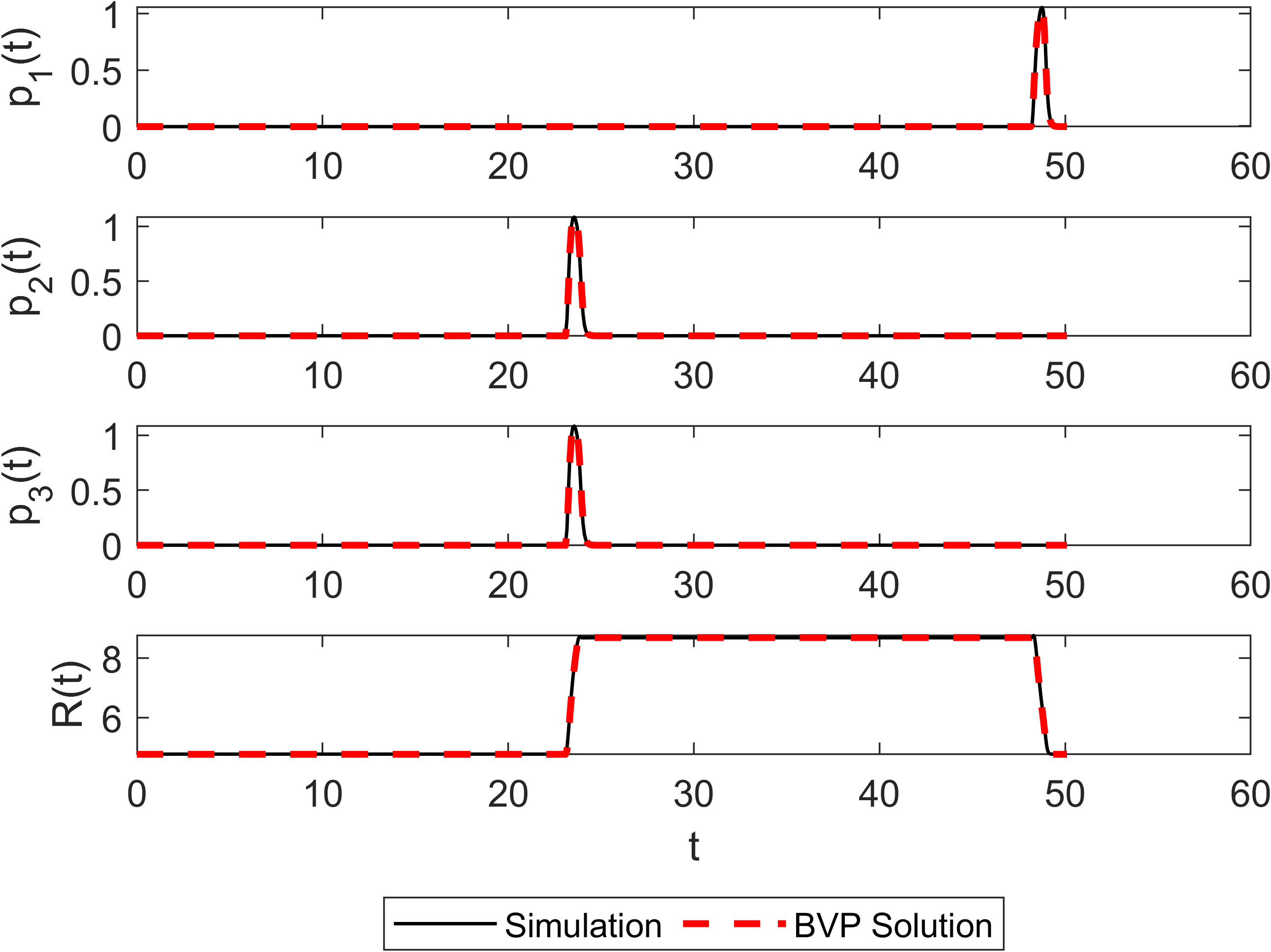}
    \caption{Boundary value problem solution for the three protein system with $\tau=25$ and $R_T=11.8$. }
    \label{fig:sametau_period25118}
\end{figure}

Next, we study the solution when $\tau=25$ at the same resource $R_T=50$. This point corresponds to a region in the response feature diagram in Fig.~\ref{fig:3ptn_amplitude} with the same amplitude as the solution in Fig.~\ref{fig:sametau_period57}. Plotting the response for these parameters between 15,000 and 15,025 time units yielded Fig.~\ref{fig:sametau_period2550}. We see that as the total resource has decreased, the periodic solution is at zero protein production rates for most of the period. Again, the protein production rates appear to oscillate in-phase for these parameters.
Lastly, we compute a periodic solution in the thin region of larger amplitude feature at low total resource. To do this, we chose $\tau=25$ and $R_T=11.8$. The periodic solution was obtained using simulation data between 15,000 and 15,025 time units and the resulting solution is shown in Fig.~\ref{fig:sametau_period25118}. It is clear that the periodic solution at these parameters is much different from the others. We see that the peak for $p_1(t)$ has shifted out of phase with the other proteins. This behavior could indicate that at extremely low resource, the best way to share that resource is to separate production of different proteins to different times, as this results in a more efficient use of resource for the cell. 

\subsection{Three Protein Summary}

The three protein system was found to behave similarly to the single protein system. Three equilibrium points were found numerically, and a subcritical Hopf bifurcation was found in the $\tau-R_T$ parameter space where all three proteins exhibited the same production time $\tau$. A large region of periodic solutions was found by numerical simulation by way of a resource limiting history function which aligns with experimental observations \cite{Tu2005}. It was found that for small delay and large total resource, the cell can produce all proteins at a constant rate at the top equilibrium point.  Due to the subcritical Hopf bifurcation of the top equilibrium point, an oscillation region appears in the parameter space which facilitates the transition from constant production to zero production just as was observed in the single protein system. Below the Hopf bifurcation curve, there is a bistability between the trivial equilibrium and the limit cycle, but slightly above the Hopf bifurcation curve there is a bistability between the top and trivial equilibrium points and the limit cycle. The trivial equilibrium was found to have a small basin of attraction and is only approached for small initial protein production rates $p_0$.

%% file: Sections/conclusion.tex
\section{Conclusion}\label{sec:conclusion}

A time delay framework was introduced for modeling protein synthesis in yeast cells. The single protein and three protein variants of this model were studied to locate regions in the parameter space where the metabolic activity contained oscillations. Three numerical methods were utilized for locating limit cycles in these systems i.e. spectral element approach applied to studying the linear stability of the systems, brute force response feature diagrams, and boundary value computation of periodic solutions to nonlinear delay equations using the spectral element method. Combining results from these three methods, the stability of each system was characterized in the parameter space where three distinct behaviors were observed. When the cell is starved of resources, it is not able to sustain protein synthesis so the protein production rates approach zero. For sufficient resources, the cell can produce at a constant rate so the system is fixed point stable, and for certain parameters between these two regions we observe limit cycle oscillations in the system. It was also observed that for the three protein model when the resources were shared and each protein had an equal production time, certain parameters at low total resource resulted in a temporal shift in the protein production rate peaks. Our simulation results are consistent with what has been observed in experiment, and our model helps argue that the observed temporal shift is a more efficient use of resources for the cell. 